\documentclass{proc-l}
\usepackage{color}
\usepackage{mathrsfs}
\usepackage{amssymb}
\usepackage{amsmath}
\usepackage{setspace}
\usepackage{makecell}
\usepackage{multirow}
\usepackage{diagbox}
\usepackage{graphicx}
\usepackage{booktabs} 
\usepackage{subcaption} 
\usepackage{float}

\newtheorem{theorem}{Theorem}[section]
\newtheorem{lemma}[theorem]{Lemma}
\newtheorem{proposition}[theorem]{Proposition}
\newtheorem{corollary}[theorem]{Corollary}

\newtheorem{question}[theorem]{Question}

\theoremstyle{definition}
\newtheorem{definition}[theorem]{Definition}
\newtheorem{example}[theorem]{Example}

\theoremstyle{remark}
\newtheorem{remark}[theorem]{Remark}

\numberwithin{equation}{section}

\begin{document}

\title[The Assouad dimension of graphs]{The Assouad Dimension of Graphs: An Affirmative Answer to Fraser's Question}



\author{Binyan Yu}
\address{\noindent{\bf Binyan Yu:} School of Mathematics and Statistics, Nanjing University of Science and Technology, 210094, Nanjing, China}
\email{yubinyan\_njust@163.com}
\email{yubinyan@njust.edu.cn}


\author{Bilel Selmi}
\address{\noindent{\bf Bilel Selmi:} Analysis, Probability and Fractals Laboratory LR18ES17, Department of Mathematics, Faculty of Sciences of Monastir,
University of Monastir, 5000, Monastir, Tunisia}
\email{bilel.selmi@fsm.rnu.tn}
\email{bilel.selmi@isetgb.rnu.tn}

\author{Jia Liu}
\address{\noindent{\bf Jia Liu:} Institute of Statistics and Applied Mathematics, Anhui University of Finance and Economics, Bengbu, 233030, China}
\email{liujia860319@163.com}

\author{Yongshun Liang}
\address{\noindent{\bf Yongshun Liang:} School of Mathematics and Statistics, Nanjing University of Science and Technology, 210094, Nanjing, China}
\email{yongshun\_liang@163.com}


\date{}

\dedicatory{}


\begin{abstract}
	This article primarily investigates the Assouad dimension of graphs of continuous functions. Initially, the Assouad spectrum of the graph of a H\"older and count-H\"older continuous function has been derived. It has been shown that the Assouad dimension of the graph of a continuous function $f:[a,b]\to \mathbb{R}$ with the count-H\"older condition is at least 
	\(
	1 + \sup \left\{ \overline{\dim}_{\mathrm{B}}L_f(y) : y \in \mathscr{R}_f[a,b]\right\}
	\)
	where $\overline{\dim}_{\mathrm{B}}$ denotes the upper box dimension, $L_f(y)$ is the level set of $f$ at level $y$, and $\mathscr{R}_f[a,b]$ is the range of $f$ on $[a,b]$. As a consequence, non-trivial lower bounds are obtained for the Assouad dimension of the graphs of classical functions such as the Weierstrass and Takagi functions. We have also employed a unified projection solely in the direction of the $y$-axis to establish lower and upper bounds for the Assouad dimension. On this basis, we have shown that the lower bound is sharp in terms of level sets, namely the upper box dimension cannot be replaced by the Assouad dimension. This provides a partial positive answer to Question 17.11.1 posed by Fraser in~\cite{Fraser2021}, and to Question 2.2 posed by Yu in~\cite{Yu2020}.
	Building on this result, the analysis has been further extended to a broader class of Weierstrass-type functions and some self-affine functions like fractal interpolation functions. Moreover, we have discussed certain circumstances under which the Assouad dimension can reach 2 with several examples constructed. A supplementary explanation of why Brownian motion has the full Assouad dimension has also been provided based on the preceding work of Howroyd and Yu \cite{Howroyd2019}.
\end{abstract}
\maketitle

\tableofcontents


\section{Introduction and main results}\label{sec1}
A significant area of investigation concerns the relationship between the analytic regularity of a continuous function \( f \colon [a, b] \to \mathbb{R} \) and the geometric regularity exhibited by its graph:
\[
\operatorname{Gr}(f, [a, b]) := \Big\{(t, f(t)) \in \mathbb{R}^2 : t \in [a, b]\Big\}.
\]
In particular, one often examines how assumptions regarding the smoothness or regularity of \( f \) impact various notions of fractal complexity associated with \( \operatorname{Gr}(f, [a, b]) \), such as different metric dimensions. 

A classical and illustrative example is the Weierstrass function \( W_{\lambda,\beta}:\mathbb{R} \to \mathbb{R} \) given by
\begin{equation}\label{eq1}
	W_{\lambda,\beta}(x) = \sum_{n=0}^{\infty} \lambda^n \cos(2\pi \beta^n x),
\end{equation}
where \( \beta > 1 \) and \( \lambda \in ( \frac{1}{\beta}, 1 ) \).
This function has been a subject of ongoing interest in the fractal geometry community, which is continuous everywhere but differentiable nowhere. This counterintuitive behavior challenges conventional concepts of differentiability and regularity in analysis. Historically important, the Weierstrass function was among the first known examples to exhibit self-affine and scale-invariant irregularity. Its graph, exhibiting intricate details at every scale, serves as a model for exploring geometric complexity. Notably, its graph has inspired extensive research aimed at determining various fractal dimensions such as the Hausdorff, packing and box dimensions. While the upper box dimension of \( \operatorname{Gr}(W_{\lambda,\beta},[a,b])\) has been investigated since the early 20\textsuperscript{th} century~\cite{Hardy1916}, the precise Hausdorff dimension was only rigorously established in 2018 by Shen~\cite{Shen2018}, who showed that \( \operatorname{Gr}(W_{\lambda,\beta},[a,b])\) has the identical Hausdorff and box dimensions
$ 2+\frac{\log\lambda}{\log\beta}$
provided $\beta$ is an integer no less than 2. More literature on its research course can be found in \cite{Besicovitch1937,Hunt1998,Hunt1991,Ren2021}.

Another fundamental example of continuous but nowhere differentiable functions is the Takagi function, which also has emerged as a significant object in fractal geometry due to its inherent self-similarity and irregularity. It serves as a canonical model for constructing and analyzing fractal sets and curves, offering valuable insights into dimension theory and the complex geometry of function graphs. Owing to its recursive structure and scaling behavior, the Takagi function also plays a key role in the study of multifractal spectra, singular measures, and functions with prescribed local oscillations, rendering it central to various aspects of the modern fractal analysis. The Takagi function is defined by
\begin{equation}\label{eq2}
	T_{\lambda, \beta}(x) = \sum_{n=0}^{\infty} \lambda^n \, \mathrm{dist}\left(\beta^n x,\mathbb{Z}\right),
\end{equation}
where \( \beta > 1 \), \( \lambda \in [ \frac{1}{\beta}, 1 ) \), and $\mathrm{dist}\left(t,\mathbb{Z}\right)$ denotes the distance from $t\in \mathbb{R}$ to the nearest integer, which is the tent map with period 1 and defined on [0,1] as
\[
\mathrm{dist}\left(t,\mathbb{Z}\right)= 
\begin{cases}
	t, &  t \in \left[0, \frac{1}{2}\right], \\[6pt]
	1 - t, &  t \in \left[\frac{1}{2}, 1\right].
\end{cases}
\]
Historically, the term “Takagi function” referred to the specific case \( T_{\frac{1}{2}, 2} \), but it is now commonly used to refer more generally to functions of the form \( T_{\lambda, \beta} \).
There has been significant interest in investigating the Hausdorff and box dimensions of the graphs of such functions. For the box dimension, it is known from \cite[Theorem 2.4]{Baranski2015} and \cite[Section 2]{Kaplan1984} that the upper box dimension of the graph of \( T_{\lambda, \beta} \) is $ 2+\frac{\log\lambda}{\log\beta}$.  
In contrast, determining the Hausdorff dimension of the graph is more challenging. For recent developments and open problems in this area, we refer to \cite{Baranski2014, Shen2018} and the references therein.

\begin{figure}[htbp] 
	\centering 
	\begin{subfigure}[b]{0.48\textwidth} 
		\centering
		\includegraphics[width=3in, height=2in]{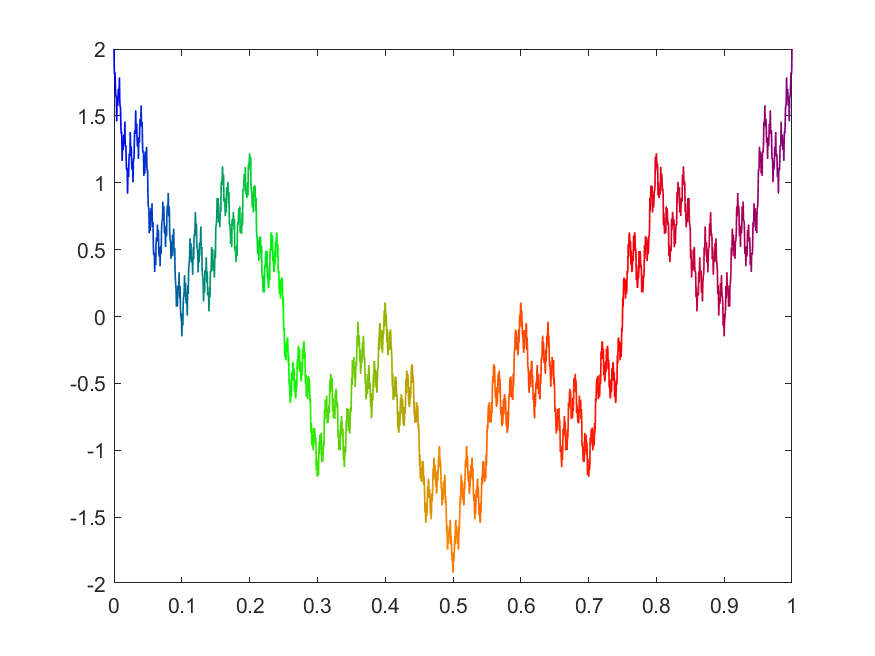} 
		\caption{The graph of $W_{\frac{1}{2},5}(x)$ on $[0,1]$.} 
		\label{fig:sub1}
	\end{subfigure}
	\hfill 
	\begin{subfigure}[b]{0.48\textwidth} 
		\centering
		\includegraphics[width=3in, height=2in]{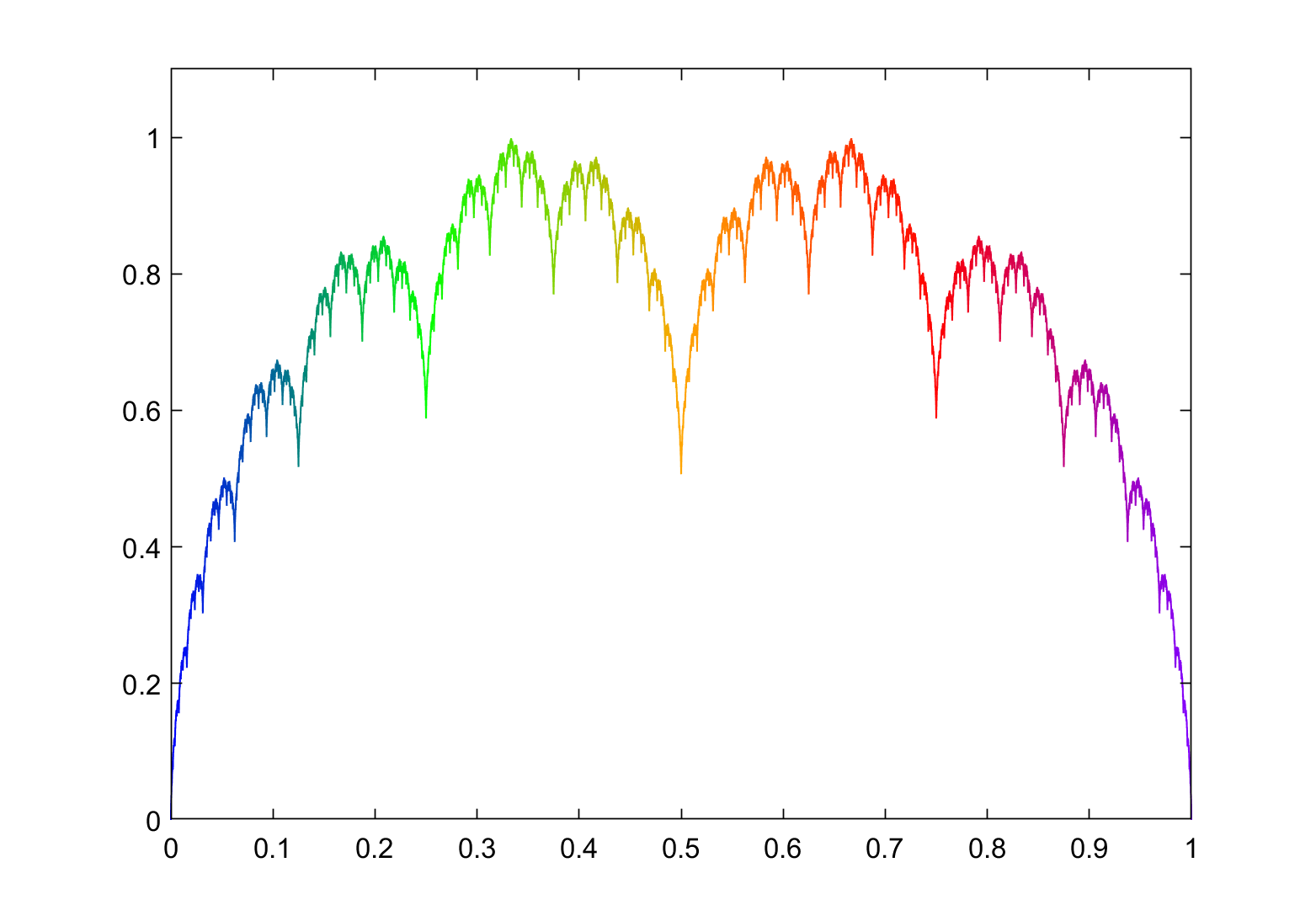} 
		\caption{The graph of $T_{\frac{2}{3},2}(x)$ on $[0,1]$.} 
		\label{fig:sub2}
	\end{subfigure}
	\caption{Graphs of the Weierstrass and Takagi functions.}\label{fig1} 
	\label{fig:parallel}
\end{figure}

In Figure \ref{fig1}, (A) and (B) display the graphs of the Weierstrass and Takagi functions on $[0,1]$ with some suitable parameters $\lambda$ and $\beta$.

This work focuses on two advanced notions of fractal dimensions: the \emph{Assouad dimension}, denoted \( \dim_{\mathrm{A}} \), introduced by Assouad in~\cite{Assouad1983}, and the \emph{Assouad spectrum}, denoted \( \dim_{\mathrm{A}}^\theta \), developed by Fraser and Yu in~\cite{FraserYu2018}. Both concepts describe the scaling behavior of coverings across a space, with the spectrum refining the Assouad dimension by incorporating a parameter \( \theta \in (0,1) \) that reflects the relationship between fine and coarse scales. In recent years, an increasing number of scholars have been attracted to the research on the Assouad dimension of certain special sets, such as popcorn graphs \cite{Chen20222}, self-affine sets \cite{FraserYu2024}, Kakeya sets \cite{Wang2025} and so on.

One outstanding problem, posed by Fraser as Question~17.11.1 in~\cite{Fraser2021}, concerns the exact value of the Assouad dimension for the graphs of classical functions such as the Weierstrass and Takagi functions. Up to now, this question still remains open. To our knowledge, the only partial progress regarding the Weierstrass function is due to Chrontsios-Garitsis and Tyson~\cite{Chrontsios-Garitsis}, who established an upper bound for its Assouad spectrum. However, for the Takagi-type function, several authors have contributed estimates and partial results on the Assouad dimension; see, for instance, Yu~\cite{Yu2020}, Anttila, Bárány, and Käenmäki~\cite{Anttila20224}, as well as the work of Chrontsios-Garitsis and Tyson~\cite{Chrontsios-Garitsis}.

Motivated by these developments, we analyze how certain regularity conditions—specifically the Hölder and count-Hölder continuity (introduced later in Section \ref{subsec2.3})—affect the lower and upper bounds for the Assouad dimension and the Assouad spectrum of \( \operatorname{Gr}(f, [a, b]) \) for a continuous function $f$ on any closed interval \( [a, b] \subset \mathbb{R} \). Our findings provide a constructive response to the open questions raised by Fraser in~\cite{Fraser2021} and by Yu in~\cite{Yu2020}, and the principal contributions of this paper are encapsulated in the upcoming main theorems.

For any continuous function $f:[a,b]\to \mathbb{R}$, our first main result establishes the relationship between the Assouad spectrum of its graph on $[a,b]$, the order $\alpha$ corresponding to its H\"older or count-H\"older continuity, and the upper box dimension (denoted $\overline{\dim}_{\mathrm{B}}$) of its level set
$$L_f(y):=\left\lbrace x\in [a,b] : f(x)=y\right\rbrace$$ at level $y\in  \mathscr{R}_f[a,b]$ where $\mathscr{R}_f[a,b]$ signifies the range of $f$ on $[a,b]$, presented below.

\begin{theorem}\label{the1.5}
	Let $ f $ be a continuous function on $[a,b]$ and $ \alpha\in(0,1) $. 
	\begin{enumerate}
		\item[(1)]   If $ f $ is $ \alpha $-order H\"older continuous on $ [a,b]$, then $$  \dim_\mathrm{A}^{\theta}\mathrm{Gr}(f,[a,b])\le\begin{cases}
			1+\displaystyle\frac{1-\alpha}{1-\theta}, \quad & \theta\in(0,\alpha),\\ \\
			2,\quad & \theta\in[\alpha,1).
		\end{cases}  $$
		\item[(2)]  If $f$ is $ \alpha $-order count-H\"older continuous on $[a,b]$, then
		$$  \dim_{\mathrm{A}}^{\theta}\mathrm{Gr}(f,[a,b])\ge \begin{cases}
			\sup\limits_{y\in \mathscr{R}_f[a,b]} \left\lbrace \overline{\dim}_{\mathrm{B}}L_f(y)\right\rbrace+\dfrac{1-\alpha}{1-\theta},\quad & \theta\in(0,\alpha),\\
			\\
			\sup\limits_{y\in \mathscr{R}_f[a,b]} \left\lbrace \overline{\dim}_{\mathrm{B}}L_f(y)\right\rbrace+1,\quad &\theta\in[\alpha,1).
		\end{cases}$$
	\end{enumerate}
\end{theorem}

\begin{figure}[htbp]
	\centering
	\begin{subfigure}[b]{0.48\textwidth}
		\centering
		\includegraphics[width=\linewidth]{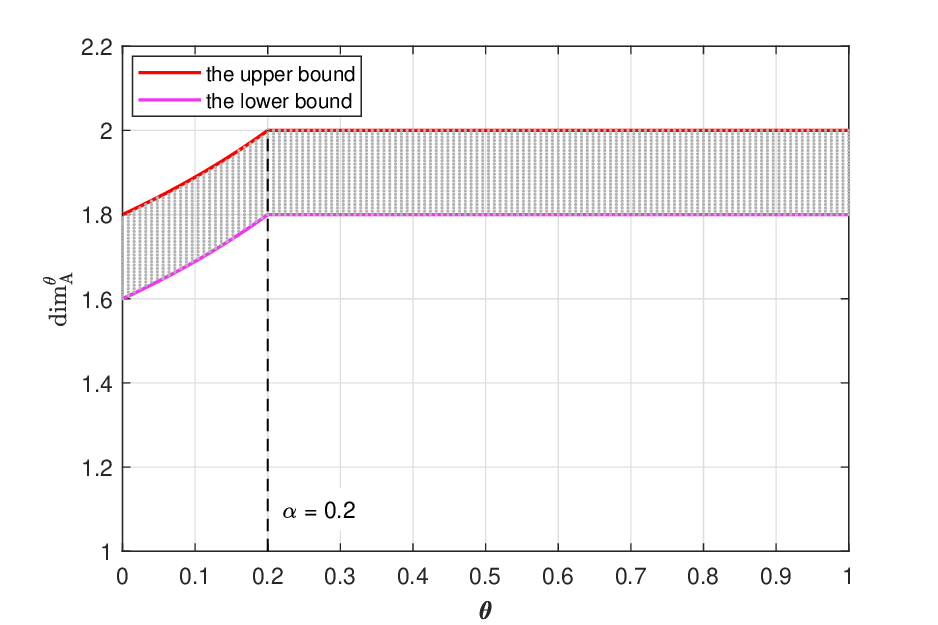}
		\caption{$\alpha=0.2$}
		\label{fig:sub1}
	\end{subfigure}
	\hfill
	\begin{subfigure}[b]{0.48\textwidth}
		\centering
		\includegraphics[width=\linewidth]{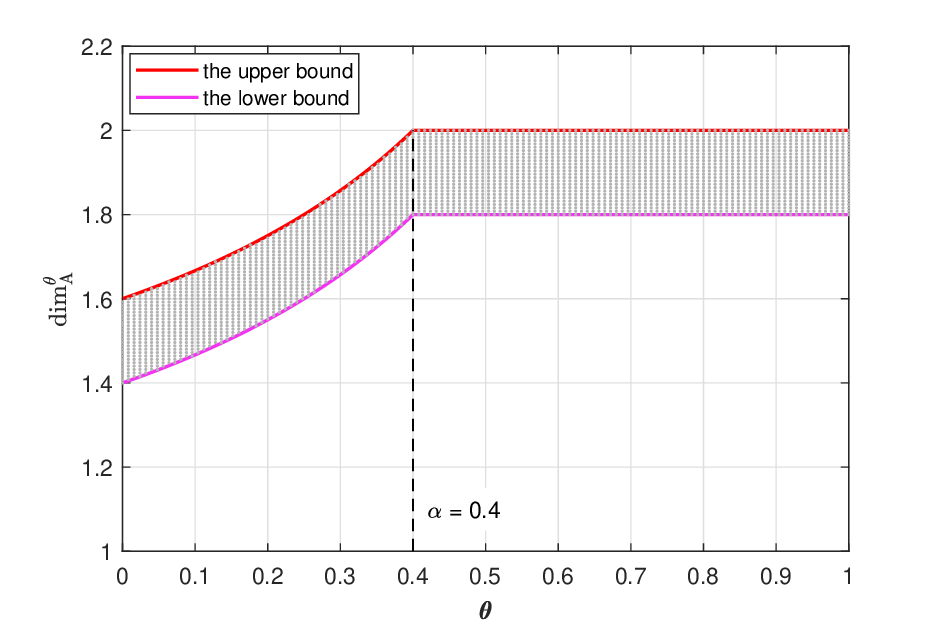}
		\caption{$\alpha=0.4$}
		\label{fig:sub2}
	\end{subfigure}
	
	\begin{subfigure}[b]{0.48\textwidth}
		\centering
		\includegraphics[width=\linewidth]{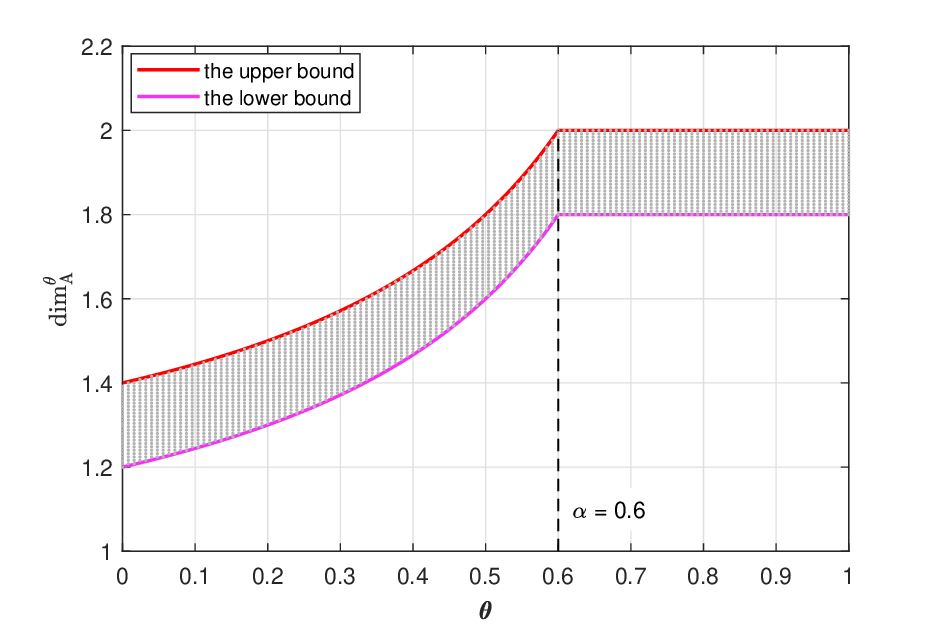}
		\caption{$\alpha=0.6$}
		\label{fig:sub3}
	\end{subfigure}
	\hfill
	\begin{subfigure}[b]{0.48\textwidth}
		\centering
		\includegraphics[width=\linewidth]{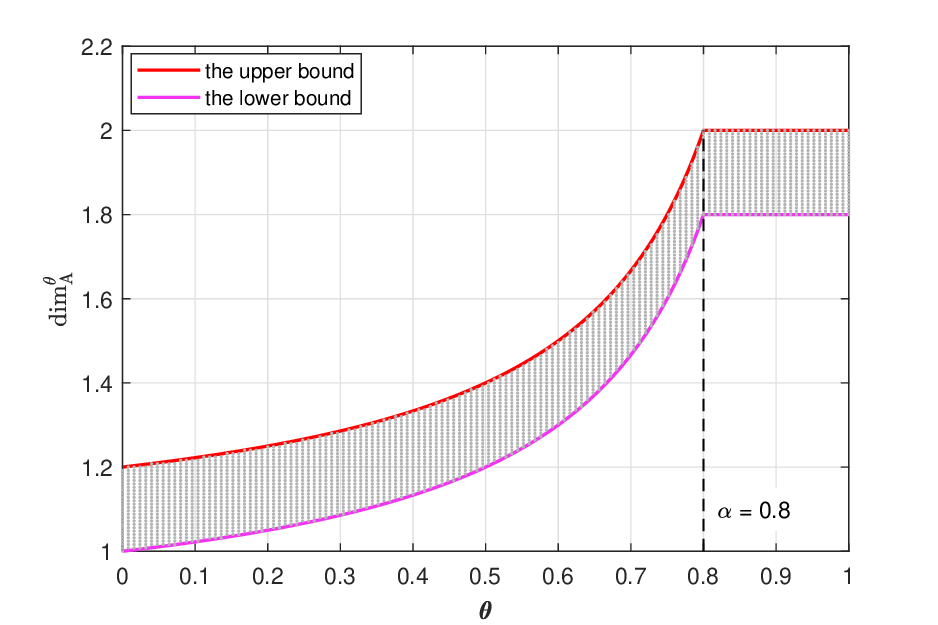}
		\caption{$\alpha=0.8$}
		\label{fig:sub4}
	\end{subfigure}
	\caption{The Assouad spectrum of the graph of any $ \alpha $-order H\"older and count-H\"older continuous function for some different values of $\alpha$.}
	\label{fig2}
\end{figure}

For the sake of intuitiveness, the Assouad spectrum of the graph of any $ \alpha $-order H\"older and count-H\"older continuous function has been visually shown in Figure \ref{fig2} (A)--(D) where we assume the supremum of the upper box dimensions of all the level sets is 0.8, and $\alpha$ is taken as 0.2, 0.4, 0.6 and 0.8, respectively. Here, the upper line represents the upper bound of the spectrum, the lower line stands for the lower bound of the spectrum, and the shaded area between these two lines means the possible range of the spectrum. Observe that the gap between the upper and lower bounds is exactly $1-\sup_{y\in \mathscr{R}_f[a,b]} \left\lbrace \overline{\dim}_{\mathrm{B}}L_f(y)\right\rbrace$ (it is 0.2 under the above assumption).

The next main result follows directly from Theorem \ref{the1.5}, which estimates the bound for the Assouad dimension of the graph of a function with the H\"older and count-H\"older continuity.

\begin{theorem}\label{the1.7}
	If $f$ is an $ \alpha $-order H\"older and count-H\"older continuous function on $[a,b]$ for $\alpha\in (0,1)$, then
	$$ \sup\limits_{y\in \mathscr{R}_f[a,b]} \left\lbrace \overline{\dim}_{\mathrm{B}}L_f(y)\right\rbrace+1\le\dim_{\mathrm{qA}}\mathrm{Gr}(f,[a,b])\le\dim_{\mathrm{A}}\mathrm{Gr}(f,[a,b])\le   2.$$
\end{theorem}

We emphasize that Theorem \ref{the1.7} still remains valid under the weaker assumption of count-H\"older continuity alone because $\dim_{\mathrm{A}}\mathrm{Gr}(f,[a,b])\le 2$ is trivial for any continuous function $f:[a,b]\to \mathbb{R}$. In addition, Theorem \ref{the1.7} seems to provide a relatively optimal lower bound, the supremum of the upper box dimensions of all level sets plus 1, for the Assouad dimension of graphs of count-H\"older continuous functions. However, it might become harder to acquire the same upper bound. 

In the present paper, we denote by $(u_{n})$ the sequence of real numbers $u_{n}$ (if necessary to stress that $n$ ranges over a specific set such as the set of positive integers $\mathbb{N}_+$, we may also denote it as $(u_{n})_{\mathbb{N}_+}$). For $ \alpha\in(0,1) $, we now define a set of binary sequences by
$$ \Lambda_{\alpha}:=\left\lbrace (r_{n},R_{n})\ \bigg|\ 0<r_n<R_n\le r_n^{\alpha}<1,\  r_{n}\searrow0,\ R_{n}\searrow0 \ \ \mathrm{and} \ \ \frac{R_{n}}{r_{n}}\nearrow\infty \ \mathrm{as} \ n\to\infty \right\rbrace.  $$
The undermentioned theorem explicitly reveals that the lower bound given in Theorem \ref{the1.7} is sharp in terms of level sets. More precisely, the upper box dimension cannot be replaced by the Assouad dimension unless we add some strict constraints to level sets.

\begin{theorem}\label{the1.3}
	Let $f$ be an $\alpha$-order count-H\"older continuous function on $[a,b]$ for $\alpha\in(0,1)$. If the following set
	$$  \begin{aligned}
		\mathcal{Y}:=\Bigg\{ y\in \mathscr{R}_{f}[a,b]\ \Bigg|\ & \mathrm{there}\ \mathrm{is} \ \mathrm{a}\ \mathrm{binary} \ \mathrm{sequence}\ (r_{n},R_{n})\in \Lambda_{\alpha}\
		\mathrm{such}\ \mathrm{that}, \\ &\hspace{1cm} \mathrm{for}\ \mathrm{all}\  n\in \mathbb{N},\
		\sup_{x\in L_{f}(y)} \#\left\lbrace  i\in \mathbb{N} \ \Big| \ \Delta_{x,n,i}\cap L_{f}(y)\neq \emptyset \right\rbrace\ge \left( \frac{R_{n}}{r_{n}}\right)^{\dim_{\mathrm{A}}L_{f}(y)} \Bigg\}
	\end{aligned} $$
	is non-empty, then 
	$$ \dim_{\mathrm{A}}\mathrm{Gr}(f,[a,b])\ge \sup_{y\in\mathcal{Y}}\left\lbrace \dim_{\mathrm{A}}L_{f}(y)\right\rbrace+1.  $$
	Here $ \Delta_{x,n,i}=\left[x-\frac{R_{n}}{2}+(i-1)r_{n},x-\frac{R_{n}}{2}+ir_{n}\right]\cap\left[x-\frac{R_{n}}{2},x+\frac{R_{n}}{2}\right] $.
\end{theorem}

Very recently, one of the most impressive findings for us is a general slicing theorem for weak tangents, which is an immediate consequence of Furstenberg’s dimension conservation result \cite[Theorem 1.6]{Furstenberg2008} and has been presented by Anttila and Rutar in \cite[Proposition A]{Anttila2024} where they have also given a short elementary proof. Here we restate this result for the completeness of our discussion.

\begin{proposition}[\cite{Anttila2024}]\label{pro1.1}
	For any non-empty compact set $ F\subset \mathbb{R}^{2} $ and any $ W\in \mathbb{RP}^{1} $, there exists an $ E\in \mathrm{Tan}(F) $ and $ x\in \mathrm{Proj}_{W}(E) $ such that
	$$ \dim_{\mathrm{H}}\left(\mathrm{Proj}^{-1}_{W}(x)\cap E \right)\ge \max\left\lbrace \dim_{\mathrm{A}}F-\dim_{\mathrm{A}}\mathrm{Proj}_{W}(F),0\right\rbrace,  $$
	where $ \mathrm{Tan}(F) $ denotes the set of all weak tangents of $ F $, $ \mathbb{RP}^{1} $ denotes the real projective space of one-dimensional subspaces of $\mathbb{R}^{2}  $ and $ \mathrm{Proj}_{W}:\mathbb{R}^{2}\to W  $ denotes the orthogonal projection onto $ W $.
\end{proposition}

Proposition \ref{pro1.1} explicitly offers an upper bound for the Assouad dimension of any non-empty compact set $ F\subset \mathbb{R}^{2} $ including graphs of continuous functions. That is, for all directions $ W $ there is a weak tangent $ E $ that makes $ \dim_{\mathrm{H}}\left(\mathrm{Proj}^{-1}_{W}(x)\cap E \right)+\dim_{\mathrm{A}}\mathrm{Proj}_{W}(F) $ the upper bound for the Assouad dimension of $ F $. For the convenience of notations, let 
$$ \mathbb{S}_{\mathrm{UB}}\left(\Omega,W\right) :=\sup\limits_{x\in\mathrm{Proj}_{W}(\Omega)} \left\lbrace \overline{\dim}_{\mathrm{B}}\left( \mathrm{Proj}_{W}^{-1}(x)\cap\Omega\right) \right\rbrace $$ 
and
$$ \mathbb{S}_{\mathrm{H}}\left(\Omega,W\right) :=\sup\limits_{x\in\mathrm{Proj}_{W}(\Omega)} \left\lbrace {\dim}_{\mathrm{H}}\left( \mathrm{Proj}_{W}^{-1}(x)\cap\Omega\right) \right\rbrace $$
for $ W\in \mathbb{RP}^{1} $ and $ \Omega\in\mathbb{R}^{2} $. Besides, we denote the $ x $-axis and the $ y $-axis, respectively, by $ \mathbb{X} $ and $ \mathbb{Y} $. Combining Theorem \ref{the1.7} and Proposition \ref{pro1.1}, we can now present the following result, which is in a unified projection form only towards the $ y $-axis direction and for both comparatively optimum lower and upper bounds.

\begin{theorem}\label{the1.9}
	If $f$ is an $ \alpha $-order count-H\"older continuous function on $[a,b]$ for $\alpha\in (0,1)$, then
	$$ \begin{aligned}
		\mathbb{S}_{\mathrm{UB}}\left(\mathrm{Gr}(f,[a,b]),\mathbb{Y}\right)+1\le\dim_{\mathrm{A}}\mathrm{Gr}(f,[a,b])\le \sup\limits_{\Omega\in\mathrm{Tan}(\mathrm{Gr}(f,[a,b]))} \mathbb{S}_{\mathrm{H}}\left(\Omega,\mathbb{Y}\right)+1.
	\end{aligned}$$
\end{theorem}

\begin{remark}
    The foregoing theorem implies another fact that the following inequality
$$ \sup\limits_{y\in \mathscr{R}_f[a,b]} \left\lbrace \overline{\dim}_{\mathrm{B}}L_f(y)\right\rbrace\le \sup\limits_{\Omega\in\mathrm{Tan}(\mathrm{Gr}(f,[a,b]))} \sup\limits_{x\in\mathrm{Proj}_{\mathbb{Y}}(\Omega)} \left\lbrace {\dim}_{\mathrm{H}}\left( \mathrm{Proj}_{\mathbb{Y}}^{-1}(x)\cap\Omega\right) \right\rbrace $$ 
holds for any function $f$ with the count-H\"older continuity. In particular, if the above formula becomes an equation, the Assouad dimension will be identical to an exact value, that is,
$$\dim_{\mathrm{A}}\mathrm{Gr}(f,[a,b])= \sup\limits_{y\in \mathscr{R}_f[a,b]} \left\lbrace \overline{\dim}_{\mathrm{B}}L_f(y)\right\rbrace+1. $$
\end{remark}

In fact, both Theorems \ref{the1.5} and \ref{the1.9} are applicable to all H\"older and count-H\"older continuous functions, and the Weierstrass and Takagi functions just serve as two representative examples of such functions whose H\"older and count-H\"older continuity will be demonstrated later in Section \ref{sec5}. Therefore, one can arrive at the following conclusion for the Weierstrass and Takagi functions.

\begin{proposition}\label{cor1.3}
	Let $\beta>\frac{1}{\lambda}>1$. Provided $\beta$ is large enough, we have the following statements.
	\begin{enumerate}
		\item[(1)] If $\theta\in (0,-\frac{\log\lambda}{\log\beta} )$, then
		$$ \mathbb{S}_{\mathrm{UB}}\left(\mathrm{Gr}(\mathcal{F}_{\lambda,\beta},[a,b]),\mathbb{Y}\right)+\dfrac{1+\frac{\log\lambda}{\log\beta}}{1-\theta}\le\dim_\mathrm{A}^{\theta}\mathrm{Gr}(\mathcal{F}_{\lambda,\beta},[a,b])\le1+\dfrac{1+\frac{\log\lambda}{\log\beta}}{1-\theta}. $$
		\item[(2)] If $\theta\in [ -\frac{\log\lambda}{\log\beta},1)$, then
		$$ \begin{aligned}
			\mathbb{S}_{\mathrm{UB}}\left(\mathrm{Gr}(\mathcal{F}_{\lambda,\beta},[a,b]),\mathbb{Y}\right)+1\le\dim_\mathrm{A}^{\theta}\mathrm{Gr}(\mathcal{F}_{\lambda,\beta},[a,b])\le\dim_{\mathrm{A}}\mathrm{Gr}(\mathcal{F}_{\lambda,\beta},[a,b])
			\le\sup\limits_{\Omega\in\mathrm{Tan}(\mathrm{Gr}(\mathcal{F}_{\lambda,\beta},[a,b]))} \mathbb{S}_{\mathrm{H}}\left(\Omega,\mathbb{Y}\right)+1.
		\end{aligned}$$
	\end{enumerate}
	Here $\mathcal{F}_{\lambda,\beta}$ represents any one of the Weierstrass function $W_{\lambda,\beta}$ defined in (\ref{eq1}) and the Takagi function $T_{\lambda,\beta}$ defined in (\ref{eq2}).
\end{proposition}

\begin{remark}
	It should be pointed out that one can take $\beta=\frac{1}{\lambda}$ in the definition of the Takagi function $T_{\lambda,\beta}$ but Proposition \ref{cor1.3} is not applicable to this situation anymore owing to the requirement for its count-H\"older continuity (we will show that it needs $\beta>\frac{1}{\lambda}$ later in Section \ref{sec5}).
\end{remark}

Additionally, several related works have appeared recently which are complementary to our results. Jiang \cite{Jiang2025} studied the Assouad dimension of the graphs of generalized Takagi functions of the form $f_{\mathbf{c},b}(x)=\sum c_n \phi(b^n x)$, and proved that its Assouad dimension equals 1 whenever $\overline{\lim}_{n\to\infty} b^n |c_n|<\infty$. For the classical Takagi functions $T_{a,b}$, this yields the same Assouad dimension 1 if and only if $a\le \frac{1}{b}$. This provides a sharp result in a special case, which is fully consistent with the general lower bounds established in our Theorem \ref{the1.7}, where the Assouad dimension is bounded below by the dimension of level sets plus one. As a particular case for the Takagi function $ T_{\lambda,2} $, Anttila et al. \cite{Anttila20224} have confirmed that $\dim_{\mathrm{A}}\mathrm{Gr}(T_{\lambda,2},[0,1])<2$ when $\frac{1}{2}<\lambda<1$. This reflects from the side that there is no level set of $ T_{\lambda,2} $ having the full upper box dimension of 1, otherwise the value of $ \dim_{\mathrm{A}}\mathrm{Gr}(T_{\lambda,2},[0,1]) $ will clock up 2 on account of the lower bound estimation in Theorem \ref{the1.7}.
Furthermore, Feng and Fraser \cite{Feng2026} investigated the Assouad spectrum and dimension of typical graphs in the little H\"older and more general little $\omega$ spaces. They showed that for a typical function in the little H\"older space of order $\alpha\in(0,1)$, the Assouad spectrum attains the value $2-\frac{\alpha-\theta}{1-\theta}$ for $\theta\le\alpha$ and $2$ for $\theta>\alpha$, and in particular the Assouad dimension is $2$. Their results demonstrate that the universal upper bounds given in our Theorem 1.1(1) are sharp in the generic sense, while our lower bounds provide a precise mechanism for this phenomenon via the dimension of level sets.

Inspired by these findings, it is natural to wonder what kind of 
 continuous functions whose graphs can achieve the full Assouad dimension of 2, and we shall deal with some particular situations for this issue in Section \ref{sec6}. At the end of this section, we put forward a relevant proposition which shows the existence of such continuous functions.

\begin{proposition}\label{pro1.3}
	Let  $ \left( u_{n}\right)_{\mathbb{N}_{+}} $ be any real number sequence satisfying the following properties:
	\begin{enumerate}
	 \item[(\romannumeral1)] $\left|u_{n+1}-u_{n}\right|\searrow 0^+ $ as $n\to\infty$;
     \item[(\romannumeral2)] $\frac{u_{n+1}-\zeta}{u_{n}-\zeta}\nearrow 1^-$ as $n\to\infty$ where $\zeta$ is the only cluster point of $ \left( u_{n}\right)_{\mathbb{N}_+}$.
\end{enumerate}
	Then the following results hold:
	\begin{enumerate}
		\item[(1)] $ \dim_{\mathrm{A}}\left\lbrace u_{n} \ \big| \ n\in \mathbb{N}_{+}\right\rbrace=1 $;
		\item[(2)] Fixing $h\in \mathbb{R}$, there exists a closed interval $[a,b]$ and a continuous function $ f:[a,b]\to \mathbb{R} $ satisfying $f(u_{n})=f(\zeta)=h$ for all $ n\in \mathbb{N}_{+}$ and $ \dim_{\mathrm{A}}\mathrm{Gr}(f,[a,b])=2 $.
		
	\end{enumerate}
\end{proposition}

\begin{remark}
	From the perspective of level sets, the statement that $f(u_{n})=f(\zeta)$ for all $ n\in \mathbb{N}_{+}$ is equivalent to the existence of a $ y_0\in \mathscr{R}_f[a,b] $ such that $$ \left\lbrace u_{n}\ \big| \ n\in \mathbb{N}_{+}\right\rbrace\cup \left\lbrace \zeta\right\rbrace\subseteq L_f(y_0).  $$
    Besides, the first claim of Proposition \ref{pro1.3} indicates that the level set $ L_f(y_0) $ has the full Assouad dimension.
\end{remark}

Proposition \ref{pro1.3} establishes the feasibility that the Assouad dimension may reach 2 if the function owns a level set which contains a certain convergent sequence. It will be verified later in Section \ref{sec6} by giving constructive examples of a special class of continuous functions with a single unbounded variation point.

The structure of the remainder of this paper is organized as follows: In the upcoming Section \ref{sec2}, we first cover basic notations and required concepts in preparation for our follow-up research. Subsequently, the proofs of our main results, Theorems \ref{the1.5}, \ref{the1.7}, \ref{the1.3} and \ref{the1.9}, have been presented in Section \ref{sec3}. Later, we extend our results to Weierstrass-type functions and self-affine functions in Section \ref{sec5}, where we also explore their H\"older and count-H\"older continuity. In Section \ref{sec6}, certain situations for the full Assouad dimension of graphs have been discussed with a class of constructive examples. Finally, we add some remarks of this work and propose several open questions for further research in Section \ref{sec7}.

\section{Preliminaries}\label{sec2}

\subsection{Notations}

We begin by introducing some basic notations which have not been illustrated in Section \ref{sec1} but will be employed throughout this paper.
\begin{enumerate}
	\item[1.] $ \mathbb{N} $ and $ \mathbb{Z} $, respectively, represent the collections of all natural numbers and integers. Let $\mathbb{N}_+=\mathbb{N}\setminus\{0\}$.
	\item[2.] For any $ n\in \mathbb{N}_+ $, $ \mathbb{R}^{n} $ stands for the $ n $-dimensional Euclidean Space. In particular, if $ n=1 $, $ \mathbb{R} $ is the collection of all real numbers. 
	\item[3.] $ \left\|\cdot\right\|  $ means the usual Euclidean distance in $ \mathbb{R}^{n} $.
	\item[4.] For any non-empty set $ \Omega\subset\mathbb{R}^{n} $, the diameter of $ \Omega $ is defined by
	$$ \left|\Omega\right|=\sup\left\lbrace \left\|t-u\right\|:t,u\in\Omega\right\rbrace.   $$
	\item[5.]   For any finite set $ \Omega\subset\mathbb{R}^{n} $, $ \# \Omega $ signifies the number of elements in $ \Omega $. 
		\item[6.] Write 
		$$  B(x,r)=\left\lbrace \zeta\in \mathbb{R}^{n} : \left\|\zeta-x\right\|\le r  \right\rbrace   $$
		as the closed ball with center $ x\in \mathbb{R}^{n} $ and radius $ r > 0 $.
        \item[7.] For $ x\in \mathbb{R} $, let $ \lceil x\rceil $ be the minimal integer no less than $x$ and $ \lfloor x\rfloor $ be the maximal integer no more than $ x $.
		\item[8.]  For any sequence $ \left( u_{n}\right)_{\mathbb{N}_+} $ and any $ \upsilon\in \mathbb{R} $, $ u_{n}\nearrow\upsilon^-  $ indicates $ \lim_{n\to\infty}u_{n}=\upsilon $ and $ u_{n}<u_{n+1} $ for all $  n\in \mathbb{N}_+ $
		while $ u_{n}\searrow \upsilon^+  $ implies $ \lim_{n\to\infty}u_{n}=\upsilon $ and $ u_{n}>u_{n+1} $ for all $  n\in \mathbb{N}_+ $.
		\item[9.] Given a function $ f: [a,b]\to \mathbb{R} $, the oscillation of $ f $ on any subinterval $ J\subseteq[a,b] $ is denoted as
		$$ \mathrm{Osc}(f,J)=\sup_{t,u\in J}\left|f(t)-f(u)\right|.  $$
		\item[10.] We denote the family of all continuous functions defined on $ [a,b] $ by $ C_{[a,b]} $.
\end{enumerate}

\subsection{Assouad dimension and Assouad spectrum}\label{subsection2.2}

In this subsection, we briefly recall basic notions of the Assouad dimension and the Assouad spectrum. The Assouad dimension, first introduced by Assouad in the 1970s \cite{Assouad1983}, reflects the ‘local characteristics’ of a set where merely covers of small neighbourhoods are considered at any one time. Let $ \Omega\subset \mathbb{R}^{n} $ be a non-empty set and let $ N_{r}(\Omega) $ be the minimal number of sets of diameter at most $ r>0 $ required to cover $ \Omega $. Then, the Assouad dimension of $ \Omega $ is defined by
$$ \begin{aligned}
	\dim_\mathrm{A} \Omega = \inf \Bigg\{ s>0 \ \Bigg|\  &\text{there exists a constant} \ C>0 \text{ such that,}\\ &\hspace{1cm} \text{ for all } 0 < r < R,\ \sup_{x \in \Omega} N_r(B(x, R) \cap \Omega) \leq C \left( \frac{R}{r} \right)^{s} \Bigg\}.
\end{aligned}$$
Building upon the definition of the Assouad dimension, the quasi-Assouad dimension proposed by L\"u and Xi \cite{Lu2016} incorporates an exponential separation between the scales
$ r $ and $ R $, which is defined as $ 	\dim_\mathrm{qA} \Omega=\lim_{\varepsilon\to0^{+}}\Psi_{\Omega}(\varepsilon) $ for $ \varepsilon\in(0,1) $ where
$$ \begin{aligned}
	\Psi_{\Omega}(\varepsilon)= \inf \Bigg\{ s>0  \ \Bigg|\ &\text{there exists a constant} \ C>0 \text{ such that,}\\ &\hspace{1cm}\text{ for all } 0 < r \le R^{1+\varepsilon}<1,\ \sup_{x \in \Omega} N_r(B(x, R) \cap \Omega) \leq C \left( \frac{R}{r} \right)^{s} \Bigg\}.
\end{aligned}$$
The Assouad spectrum introduced by Fraser and Yu \cite{FraserYu2018} is a function of the parameter $\theta\in(0,1)$,
which can be an interpolation between the upper box dimension and the quasi-Assouad dimension as the parameter varies. To be specific, define
$$ \begin{aligned}
	\dim_\mathrm{A}^{\theta} \Omega = \inf \Bigg\{ s>0 \ \Bigg|\  &\text{there exists a constant}\ C>0 \text{ such that,}\\ &\hspace{1cm} \text{ for all } 0 < r < 1,\ \sup_{x \in \Omega} N_r(B(x, r^{\theta}) \cap \Omega) \leq C \left( \frac{r^{\theta}}{r} \right)^{s} \Bigg\}
\end{aligned}$$ 
and the Assouad spectrum of $\Omega$ is a family of values $\left\lbrace 	\dim_\mathrm{A}^{\theta} \Omega : \theta\in(0,1) \right\rbrace $.

\begin{remark}
	To estimate the Assouad dimension of the graph of a continuous function $ f: [a,b]\to \mathbb{R} $, note that for $x\in[a,b]$ the covering number $ N_{r}(B((x,f(x)), R) \cap \mathrm{Gr}(f,[a,b])) $ in the above definitions can be replaced by the least number of cubes of sidelength $r$ required to cover $B((x,f(x)), R) \cap \mathrm{Gr}(f,[a,b])$ without changing the values for dimensions $\dim_\mathrm{A}\mathrm{Gr}(f,[a,b])$, $\dim_\mathrm{qA}\mathrm{Gr}(f,[a,b])$ and $\dim_\mathrm{A}^{\theta} \mathrm{Gr}(f,[a,b])$. As an equivalent definition, the ball $B((x,f(x)), R)$ can also be replaced by a cube of sidelength $R$ centred at $(x,f(x))$, which can be written as
	$$Q((x,f(x)),R)=\left\lbrace (y_1,y_2)  \ \bigg|\ \vert y_1-x\vert\le \frac{R}{2},\ \vert y_2-f(x)\vert\le \frac{R}{2} \right\rbrace$$ 
	for the convenience of notations in later parts.
\end{remark}

Next, we present certain fundamental properties of the Assouad dimension and the Assouad spectrum, which will be devoted to the subsequent proofs. For any bounded set $ \Omega\subset \mathbb{R}^{n} $, the following results hold for $\dim_\mathrm{A}^{\theta} \Omega$: 
\begin{enumerate}
	\item[(1)] The map $\theta\mapsto\dim_\mathrm{A}^{\theta} \Omega$ is continuous (even locally Lipschitz) in $\theta\in(0,1)$;
	\item[(2)] $\dim_\mathrm{A}^{\theta} \Omega\to \overline{\dim}_\mathrm{B}\Omega$ as $\theta\to0^+$;
	\item[(3)] $\dim_\mathrm{A}^{\theta} \Omega\to \dim_\mathrm{qA}\Omega$ as $\theta\to1^-$.
\end{enumerate}
Here $\overline{\dim}_\mathrm{B}\Omega$ represents the upper box dimension of $\Omega$. It should be noted that there exists an essential relationship between those dimensions, that is, for all $\theta\in(0,1)$,
\begin{equation}\label{eq3}
	\overline{\dim}_\mathrm{B}\Omega\le\dim_\mathrm{A}^{\theta} \Omega\le\dim_\mathrm{qA} \Omega \le \dim_\mathrm{A}\Omega.
\end{equation} 
In fact, for many sets of interest, it holds $\dim_\mathrm{qA} \Omega= \dim_\mathrm{A}\Omega$ so that the entire range $[\overline{\dim}_\mathrm{B}\Omega,\dim_\mathrm{A}\Omega]$ is just described by the Assouad spectrum. 
For more details and recent work on the Assouad dimension and the Assouad spectrum, the interested readers can consult, for instance, \cite{Anttila20223, Banaji2024, Fraser2014, Fraser2019, FraserYu2023}.

\subsection{H\"older and count-H\"older continuous functions}\label{subsec2.3}

Regarded as a regularity condition for continuous functions, the notion of the H\"older and count-H\"older continuity at a point has been given below.

\begin{definition}\label{def6.0}
    Let $ f\in C_{[a,b]} $ and $ \alpha\in(0,1) $. 
    \begin{enumerate}
		\item[(1)] $ f $ is said to be $ \alpha $-order H\"older continuous at $ t\in[a,b] $ if there are certain constants $ H^{+}(f)>0 $ and $ \delta_{0}>0 $ such that for all $ \delta\in (0,\delta_{0}] $,
		\begin{equation}\label{eq6.1}
		\mathrm{Osc}(f,[t-\delta,t+\delta]\cap[a,b])\le H^{+}(f)\cdot\delta^{\alpha}.
		\end{equation}
		\item[(2)] $ f $ is said to be $ \alpha $-order count-H\"older continuous at $ t\in[a,b] $ if there are certain constants $ H^{-}(f)>0 $ and $ \delta_{0}>0 $ such that for all $ \delta\in (0,\delta_{0}] $,
		\begin{equation}\label{eq6.2}
			\mathrm{Osc}(f,[t-\delta,t+\delta]\cap[a,b])\ge H^{-}(f)\cdot\delta^{\alpha}.
		\end{equation}
		\item[(3)] $ f $ is said to be $ \alpha $-order H\"older and count-H\"older continuous at $ t\in[a,b] $ if it satisfies both of conditions (\ref{eq6.1}) and (\ref{eq6.2}).
	\end{enumerate}
\end{definition}

On the basis of Definition \ref{def6.0}, we can introduce the following concept of H\"older and count-H\"older continuous functions.

\begin{definition}\label{def6.1}
	Let $ f\in C_{[a,b]} $ and $ \alpha\in(0,1) $. 
	\begin{enumerate}
		\item[(1)] $ f $ is termed an $ \alpha $-order H\"older continuous function on $ [a,b] $ if the condition (\ref{eq6.1}) holds
		for all $ t\in [a,b] $.
		\item[(2)] $ f $ is termed an $ \alpha $-order count-H\"older continuous function on $ [a,b] $ if the condition (\ref{eq6.2}) holds
		for all $ t\in [a,b] $.
		\item[(3)] $ f $ is called an $ \alpha $-order H\"older and count-H\"older continuous function on $ [a,b] $ if both conditions (\ref{eq6.1}) and (\ref{eq6.2}) hold for all $ t\in [a,b] $.
	\end{enumerate}
\end{definition}

\begin{remark}
	It is worth mentioning that the values of H\"older and count-H\"older coefficients $ H^{+}(f)$ and $H^{-}(f) $ are not unique. Actually, if we can find two positive constants $ \mathcal{H}_{1} $ and $ \mathcal{H}_{2} $ meeting the conditions (\ref{eq6.1}) and (\ref{eq6.2}), respectively, then all $ H^{+}(f)\in[\mathcal{H}_{1},+\infty) $ and $ H^{-}(f)\in(0,\mathcal{H}_{2}] $ will still keep the results correct.
\end{remark}

\begin{remark}
	According to Definition \ref{def6.1}, one can directly derive an alternative judgment for H\"older and count-H\"older continuous functions as follows:
	\begin{enumerate}
		\item[(1)] Suppose that there exists a certain constant $ H^{+}(f)>0 $ such that
		\begin{equation}\label{eq100}
			\left|f(t)-f(u)\right|\le H^{+}(f)\cdot\left|t-u\right|^{\alpha}
		\end{equation}
		for all $ t,u\in [a,b] $. Then, $ f $ is $ \alpha $-order H\"older continuous on $ [a,b] $.
		\item[(2)] Suppose that there are certain constants $ H^{-}(f)>0 $ and $ \delta_{0}>0 $ satisfying the following property: for each $ t\in [a,b] $ and $ \delta\in(0,\delta_{0}] $, there always exists a $ t_{0}\in [a,b] $ such that $ \left|t-t_{0}\right|\le\delta  $ and 
		\begin{equation}\label{eq101}
			\left|f(t)-f(t_{0})\right|\ge H^{-}(f)\cdot\delta^{\alpha}. 
		\end{equation}
		Then, $ f $ is $ \alpha $-order count-H\"older continuous on $ [a,b] $.
	\end{enumerate} 
\end{remark}

In \cite{Falconer2014}, Falconer has built the following dimensional result for H\"older and count-H\"older continuous functions.

\begin{proposition}[\cite{Falconer2014}]\label{rem1.3}
	For $ \alpha\in(0,1) $, the graph of an $ \alpha $-order H\"older and count-H\"older continuous function on $ [a,b] $ has the box dimension $ 2-\alpha $. 
\end{proposition}

We now state a fundamental conclusion for the Weierstrass and Takagi functions.

\begin{proposition}\label{pro2.6}
	For $\beta>\frac{1}{\lambda}>1$, the Weierstrass function $ W_{\lambda,\beta} $ defined in (\ref{eq1}) and the Takagi function $ T_{\lambda,\beta} $ defined in (\ref{eq2}) are both H\"older and count-H\"older continuous functions with the order $ -\frac{\log\lambda}{\log\beta} $ provided $\beta$ is large enough.
\end{proposition}

However, the calculation of H\"older and count-H\"older coefficients is not simple at all, which is the key to determining the H\"older and count-H\"older continuity. Later in Section \ref{sec5}, we will demonstrate how to acquire a pair of them through a class of H\"older and count-H\"older continuous functions and the proof of Proposition \ref{pro2.6} will be provided as well.

\section{Proofs of the main results}\label{sec3}

In this section, we prove our main results, Theorems \ref{the1.5}, \ref{the1.7}, \ref{the1.3} and \ref{the1.9}. 

\subsection{Proof of Theorem \ref{the1.5}(1)}
We shall start the proof of the first claim of Theorem \ref{the1.5} for the upper bound of the Assouad spectrum of H\"older continuous functions. To prove this result, the forthcoming basic lemma is required at first.

\begin{lemma}\label{lem5}
	Let $ f\in C_{[a,b]} $. For all $ t\in[a,b] $, we have
	$$ N_{r}\left(Q\left((t,f(t)),R\right)\cap \mathrm{Gr}(f,[a,b])\right)\le \sum_{i= 1}^{\lceil\frac{d-c}{r}\rceil}\left\lbrace 2+\frac{\min\left\lbrace R,\mathrm{Osc}\left(f,\Delta_{i}\right) \right\rbrace }{r}\right\rbrace,
	$$
	where $0 <r<R$, $ c=\max\left\lbrace t-\frac{R}{2},a \right\rbrace  $, $ d=\min\left\lbrace t+\frac{R}{2},b \right\rbrace  $, and $ \Delta_{i}=[c + (i-1)r,c+ir]\cap[c,d] $ for $ i=1,2,\cdots,\lceil\frac{d-c}{r}\rceil $.
\end{lemma}

\begin{proof}
	Since $ f $ is continuous on $ [a,b] $, the number of cubes of sidelength $ r $ in the column above the interval
	$ \Delta_{i} $ which can intersect $ Q\left((x,f(x)),R\right)\cap \mathrm{Gr}(f,[a,b]) $ is at most 
	$$ 2+\frac{\min\left\lbrace R,\mathrm{Osc}\left(f,\Delta_{i}\right) \right\rbrace }{r}. $$
	Summing over all such intervals just gives this estimation. 
\end{proof}

With the help of Lemma \ref{lem5}, we now give the following proof of Theorem \ref{the1.5}(1).

\begin{proof}[Proof of Theorem \ref{the1.5}(1)]
	Since $ f(x) $ is an $ \alpha $-order H\"older continuous function on $ [a,b] $, it is immediate from (\ref{eq6.1}) that there are certain constants $ H^{+}(f)>0 $ and $\delta_0>0$ such that
	$$ 	\begin{aligned}
		\mathrm{Osc}(f,[p,q])\le  \frac{H^{+}(f)}{2^{\alpha}}\cdot (q-p)^{\alpha}
	\end{aligned} $$
	for any closed subinterval $ [p,q]\subseteq [a,b] $ satisfying $q-p\le 2\delta_0$. For all $0 < r <\min\left\lbrace2\delta,1\right\rbrace$ and all $ x\in[a,b] $, suppose that $ [\upsilon,\omega]=\big[x-\frac{r^{\theta}}{2},x+\frac{r^{\theta}}{2}\big]\cap[a,b] $. Therefore, it follows from Lemma \ref{lem2.1} that
	\begin{align*}
		N_{r}\left(Q\left((x,f(x)),r^{\theta}\right)\cap \mathrm{Gr}(f,[a,b])\right)
		&\le \sum_{\ell=1}^{\lceil\frac{\omega-\upsilon}{r}\rceil}\left\lbrace 2+\frac{\min\left\lbrace r^{\theta},\mathrm{Osc}\left(f,[\upsilon + (\ell-1) r,\upsilon+\ell r]\cap[\upsilon,\omega]\right) \right\rbrace }{r}\right\rbrace\\
		&\le\left(1+\frac{r^{\theta}}{r}\right)\cdot\frac{\left(\frac{H^{+}(f)}{2^{\alpha}}+3\right)\cdot \min\left\lbrace r^{\theta},r^{\alpha}\right\rbrace}{r}\\
		&\le C_{1}\cdot r^{\theta-2}\cdot\min\left\lbrace r^{\theta},r^{\alpha}\right\rbrace,
	\end{align*}
	where $ C_{1}=2^{1-\alpha}H^{+}(f)+6>0 $. Thus
	$$ N_{r}\left(Q\left((x,f(x)),r^{\theta}\right)\cap \mathrm{Gr}(f,[a,b])\right)\le C_{1}\cdot r^{\theta-2+\alpha}=  C_{1}\cdot \left( \frac{r^{\theta}}{r}\right)^{1+\frac{1-\alpha}{1-\theta}} $$
	when $\theta\in(0,\alpha)$ and 
	$$ N_{r}\left(Q\left((x,f(x)),r^{\theta}\right)\cap \mathrm{Gr}(f,[a,b])\right)\le C_{1}\cdot r^{2\theta-2}=  C_{1}\cdot \left( \frac{r^{\theta}}{r}\right)^{2} $$
	when $\theta\in[\alpha,1)$, implying that the required result holds.
\end{proof}

\subsection{Proof of Theorem \ref{the1.5}(2)}
We now turn towards the proof of the second claim of Theorem \ref{the1.5}. Given a function $f\in C_{[a,b]}$, it is relatively more complicated to determine the lower bound of the Assouad spectrum of its graph due to the local restrictions by $Q\left((x,f(x)),R\right)$. To this end, we shall introduce the level sets of $f$ in calculating the lower bound of $N_{r}\left(Q\left((x,f(x)),R\right)\cap \mathrm{Gr}(f,[a,b])\right)$ and set up the undermentioned estimation corresponding to Lemma \ref{lem5}.

\begin{lemma}\label{lem2.0}
	Let $ f\in C_{[a,b]} $. For all $ t\in[a,b] $ and $ \lambda\in(0,1) $, it holds that
	$$ 	 \begin{aligned}
		&N_{r}\left(Q\left((t,f(t)),R\right)\cap \mathrm{Gr}(f,[a,b])\right)\\&\hspace{2cm}\ge \sum_{i =1}^{\lceil\frac{d-c}{r}\rceil}\mathbb{I}\left(\Delta_{i}\cap L_f\left(f(t)+\frac{2\lambda-1}{2}R\right)\right)\cdot\max\left\lbrace 1,\frac{\min\left\lbrace \frac{1-\left|2\lambda-1\right| }{2} R, \mathrm{Osc}\left(f,\Delta_{i}\right)\right\rbrace }{r}\right\rbrace,
	\end{aligned}
	$$
	where $0 <r<R$, $ c=\max\left\lbrace t-\frac{R}{2},a \right\rbrace  $, $ d=\min\left\lbrace t+\frac{R}{2},b \right\rbrace  $, $ \Delta_{i}=[c + (i-1)r,c+ir]\cap[c,d] $ for $ i=1,2,\cdots,\lceil\frac{d-c}{r}\rceil$ and $\mathbb{I}$ is the indicator function for a set $\Omega\subset\mathbb{R}^n$ defined by
	$$ \mathbb{I}\left(\Omega\right) =\begin{cases}
		1,\quad &\Omega\neq \emptyset,\\\\
		0,\quad & \Omega= \emptyset.    \end{cases}$$
\end{lemma}

\begin{proof}
	Given $ t\in[a,b] $ and $ \lambda\in(0,1) $, we observe that 
	$$ L_f\left(f(t)+\frac{2\lambda-1}{2}R\right)\cap\left[c,d\right]\subset\mathrm{Proj}_{\mathbb{X}}\left(Q\left((t,f(t)),R\right)\cap \mathrm{Gr}(f,[a,b])\right),  $$
thus $N_{r}\left(Q\left((t,f(t)),R\right)\cap \mathrm{Gr}(f,[a,b])\right)$ must be no less than the sum of cubes of sidelength $ r $ in the column above all the interval
	$ \Delta_{i}$ that can intersect $ L_f\left(f(t)+\frac{2\lambda-1}{2}R\right) $. Since $ f $ is continuous on $ [a,b] $, the number of such cubes above each required interval
	$ \Delta_{i} $ which intersects $ Q\left((t,f(t)),R\right)\cap \mathrm{Gr}(f,[a,b]) $ is at least 
	$$ \max\left\lbrace 1,\frac{\min\left\lbrace \lambda R,(1-\lambda)R,\mathrm{Osc}\left(f,\Delta_{i}\right)\right\rbrace }{r}\right\rbrace, $$
	and the number of all the required interval $ \Delta_{i} $ is at least
	$$ \sum_{i =1}^{\lceil\frac{d-c}{r}\rceil}\mathbb{I}\left(\Delta_{i}\cap L_f\left(f(t)+\frac{2\lambda-1}{2}R\right)\right).$$
	Note that $$ \min\left\lbrace \lambda,1-\lambda \right\rbrace=\frac{1-\left|2\lambda-1\right| }{2},  $$
    then summing over all such cubes just leads to this estimation.
\end{proof}

The next lemma is a corollary following from Lemma \ref{lem2.0} as long as we particularly take $ \lambda=\frac{1}{2} $.

\begin{lemma}\label{lem2.1}
	Let $ f\in C_{[a,b]} $. For all $ t\in[a,b] $, it holds
	$$ 	 \begin{aligned}
		N_{r}\left(Q\left((t,f(t)),R\right)\cap \mathrm{Gr}(f,[a,b])\right)
		\ge \sum_{i =1}^{\lceil\frac{d-c}{r}\rceil}\mathbb{I}\left(\Delta_{i}\cap L_f(f(t))\right)\cdot\max\left\lbrace 1,\frac{\min\left\lbrace \frac{R}{2},\mathrm{Osc}\left(f,\Delta_{i}\right)\right\rbrace }{r}\right\rbrace,
	\end{aligned}
	$$
	where $r$, $ R $, $ c $, $ d $ and $ \Delta_{i}$ are as shown in Lemma \ref{lem2.0}.
\end{lemma}

In the light of Lemma \ref{lem2.1}, together with the analytical property of count-H\"older continuous functions, we are able to prove Theorem \ref{the1.5}(2) as follows.

\begin{proof}[Proof of Theorem \ref{the1.5}(2)]
	For the sake of notations, write
	$$ \eta=\sup_{y\in \mathscr{R}_f[a,b]} \left\lbrace \overline{\dim}_{\mathrm{B}}L_f(y)\right\rbrace.$$
	Here we only need to consider the case when $ 0<\eta\le 1 $ since it is trivial that
	$$ \dim_{\mathrm{A}}\mathrm{Gr}(f,[a,b])\ge 1=\sup_{y\in \mathscr{R}_f[a,b]} \left\lbrace \overline{\dim}_{\mathrm{B}}L_f(y)\right\rbrace+1 $$
	for any $f\in C_{[a,b]}$ when $ \eta=0 $. Arbitrarily given an $\varepsilon>0$, there exists a $y_0\in \mathscr{R}_f[a,b]$ such that
	$$ \overline{\dim}_{\mathrm{B}}L_f(y_0)\ge \eta-\frac{\varepsilon}{2}.$$
	Let $s=\overline{\dim}_{\mathrm{B}}L_f(y_0)$. Then by (\ref{eq3}), $\dim_{\mathrm{A}}^{\theta}L_f(y_0)> s-\frac{\varepsilon}{2}$ 
    for all $\theta\in (0,1)$. Hence, given a $\theta\in (0,1)$, there is a sequence of points $x_n\in L_f(y_0)$ and a sequence of scales $0<r_n<1$ such that $r_n\searrow 0 $ as $n\to\infty$ and
	\begin{equation}\label{eq19}
		N_{r_n}\left(Q\left((x_n,y_0),r_n^{\theta}\right)\cap L_f(y_0)\right)\ge  \left( \frac{r_n^{\theta}}{r_n}\right)^{s-\frac{\varepsilon}{2}}\ge \left( \frac{r_n^{\theta}}{r_n}\right)^{\eta-\varepsilon}
	\end{equation} 
	for all $n\in \mathbb{N}$. Set $$ [c_n,d_n]=\left[x_n-\frac{r_n^{\theta}}{2},x_n+\frac{r_n^{\theta}}{2}\right]\cap [a,b]=Q\left((x_n,y_0),r_n^{\theta}\right)\cap L_f(y_0) $$ and
	$$ \Delta_{n,j}=[c_n + jr_n,c_n+(j + 1)r_n]$$
	for all $ j=1,2,\cdots,\lceil\frac{d_n-c_n}{r_n}\rceil$. Since $ f(x) $ is $ \alpha $-order count-H\"older continuous on $ [a,b] $, there are certain constants $ H^{-}(f)>0 $ and $ \delta_{0}>0 $ such that $$\mathrm{Osc}(f,[t-\delta,t+\delta]\cap[a,b])\ge H^{-}(f)\cdot\delta^{\alpha}$$
    for all $t\in[a,b]$ and $\delta\in(0,\delta_0]$.
	Moreover, there exists an $m_0\in \mathbb{N}$ such that 
	$ r_n\le 2\delta_0 $ for all $ n\ge m_0$. Thus
	\begin{equation}\label{eq10}
		\begin{aligned}
			\mathrm{Osc}\left(f,\Delta_{n,j}\right)\geq\frac{H^{-}(f)}{2^{\alpha}}\cdot r_n^{\alpha}
		\end{aligned}
	\end{equation}
	for all $ j=1,2,\cdots,\lceil\frac{d_n-c_n}{r_n}\rceil-1$ and all $ n\ge m_0$. Note that
	$$ \sum_{j =1}^{\lceil\frac{d_n-c_n}{r_n}\rceil-1}\mathbb{I}\left(\Delta_{n,j}\cap L_f(y_0)\right)\ge N_{r_n}\left([c_n,d_n]\cap L_f(y_0)\right)-1\ge \frac{1}{2}N_{r_n}\left([c_n,d_n]\cap L_f(y_0)\right).$$
	Therefore, combined with Lemma \ref{lem2.1}, (\ref{eq19}) and (\ref{eq10}), it can be estimated that for all $ n\ge m_0$,
	$$ \begin{aligned}
		N_{r_n}\left(Q\left((x_n,y_0),r_n^{\theta}\right)\cap \mathrm{Gr}(f,[a,b])\right)&	\ge \sum_{j =1}^{\lceil\frac{d_n-c_n}{r_n}\rceil - 1}\mathbb{I}\left(\Delta_{n,j}\cap L_f(y_0)\right)\cdot\frac{\min\left\lbrace \frac{r_n^{\theta}}{2},\mathrm{Osc}\left(f,\Delta_{n,j}\right)\right\rbrace }{r_n}\\
		&\ge   \frac{1}{2}\left( \frac{r_n^{\theta}}{r_n}\right)^{\eta-\varepsilon}\cdot \frac{\min\left\lbrace \frac{1}{2}r_n^{\theta},\frac{H^{-}(f)}{2^{\alpha}}\cdot r_n^{\alpha}\right\rbrace }{r_n}\\
		&= C_{2}\cdot  r_n^{(\eta-\varepsilon)(\theta-1)-1}\cdot\min\left\lbrace r_n^{\theta}, r_n^{\alpha}\right\rbrace,
	\end{aligned}$$
	where $C_2=\min\left\lbrace \frac{1}{4}, \frac{H^{-}(f)}{2^{\alpha+1}}\right\rbrace>0$. Consequently, one has
	$$ N_{r_n}\left(Q\left((x_n,y_0),r_n^{\theta}\right)\cap \mathrm{Gr}(f,[a,b])\right)\ge C_{2}\cdot  r_n^{(\eta-\varepsilon)(\theta-1)-1+\alpha}= C_{2}\cdot  \left( \frac{r_n^{\theta}}{r_n}\right)^{\eta+\frac{1-\alpha}{1-\theta}-\varepsilon}$$
	when $\theta\in(0,\alpha)$ and
	$$ N_{r_n}\left(Q\left((x_n,y_0),r_n^{\theta}\right)\cap \mathrm{Gr}(f,[a,b])\right)\ge C_{2}\cdot  r_n^{(\eta-\varepsilon)(\theta-1)-1+\theta}= C_{2}\cdot  \left( \frac{r_n^{\theta}}{r_n}\right)^{\eta+1-\varepsilon}$$
	when $\theta\in[\alpha,1)$. Because of the arbitrariness of $\varepsilon$, it can be inferred that
	$$  \dim_{\mathrm{A}}^{\theta}\mathrm{Gr}(f,[a,b])\ge \begin{cases}
		\eta+\dfrac{1-\alpha}{1-\theta},\quad & \theta\in(0,\alpha),\\
		\\
		\eta+1,\quad &\theta\in[\alpha,1),
	\end{cases}$$
	leading to our desired result.
\end{proof}

\subsection{Proof of Theorem \ref{the1.3}}
Given $f\in C_{[a,b]}$, $ y\in \mathscr{R}_{f}[a,b] $, $x\in L_f(y)$ and $0<r_{n}<R_{n}<1$, we first give the following two assumptions:
	\begin{enumerate}
	\item[$ (\mathrm{A1}) $]  There is a constant $ K_{1}>0 $ such that
	$$  \#\left\lbrace i\in \mathbb{N} \ \Big|\ \Delta_{x,n,i}\cap L_f\left(y\right)\neq\emptyset\right\rbrace\ge K_{1}\cdot \left( \frac{R_n}{r_n}\right)^{\dim_{\mathrm{A}}L_{f}\left(y\right)};     $$ 
	\item[$ (\mathrm{A2}) $] There is a constant $ K_{2}>0 $ such that  \begin{equation}\label{eq94}
	    \mathrm{Osc}\left(f,\Delta_{x,n,i}\right)\ge K_{2}\cdot R_{n}
	\end{equation}  for all $ i=1,2,\cdots,\big\lceil\frac{d_{x,n}-c_{x,n}}{r_{n}}\big\rceil$.
\end{enumerate}
Here $ c_{x,n}=\max\left\lbrace x-\frac{R_{n}}{2},a \right\rbrace  $, $ d_{x,n}=\min\left\lbrace x+\frac{R_{n}}{2},b \right\rbrace  $ and $ \Delta_{x,n,i}=[c_{x,n} + (i-1)r_{n},c_{x,n}+ir_{n}]\cap[c_{x,n},d_{x,n}] $ for $ i=1,2,\cdots,\big\lceil\frac{d_{x,n}-c_{x,n}}{r_{n}}\big\rceil$. 
Before setting up the proof of Theorem \ref{the1.3}, the next lemma is auxiliary and indispensable.

\begin{lemma}\label{lem3.4}
   Let $f\in C_{[a,b]}$. If there is a $ y\in \mathscr{R}_{f}[a,b] $, a sequence of points $x_n\in L_f(y)$ and two sequences of scales $0<r_{n}<R_{n}<1$ satisfying both assumptions $\mathrm{(A1)}$ and $\mathrm{(A2)}$ for large enough $n$, then
    $$ N_{r_{n}}\left(Q\left((x_{n},y),R_{n}\right)\cap \mathrm{Gr}(f,[a,b])\right)
		\ge C\cdot\left(\frac{R_{n}}{r_{n}}\right)^{\dim_{\mathrm{A}}L_{f}(y)+1},  $$
        where $ C $ is a certain positive constant.
\end{lemma}

\begin{proof} 
 Under this premise, there are two constants $ K_{1},K_2>0 $, a $ y\in \mathscr{R}_{f}[a,b] $, a sequence of points $x_n\in L_f(y)$ and two sequences of scales $0<r_{n}<R_{n}<1$ such that
	$$ \sum_{i=1}^{\big\lceil\frac{d_{x_n,n}-c_{x_n,n}}{r_{n}}\big\rceil}\mathbb{I}\left(\Delta_{x_n,n,i}\cap L_{f}(y)\right)=\#\left\lbrace i\in \mathbb{N} \ \Big|\ \Delta_{x_n,n,i}\cap L_f\left(y\right)\neq\emptyset\right\rbrace \ge K_{1}\cdot \left( \frac{R_n}{r_n}\right)^{\dim_{\mathrm{A}}L_{f}\left(y\right)}$$
and
\begin{equation*}
	\mathrm{Osc}\left(f,\Delta_{x_n,n,i}\right)\ge K_{2}\cdot R_{n}
\end{equation*}
for all $ i=1,2,\cdots,\big\lceil\frac{d_{x_n,n}-c_{x_n,n}}{r_{n}}\big\rceil$ and large enough $n$. Then, it follows from Lemma \ref{lem2.1} that
$$ \begin{aligned}
N_{r_n}\left(Q\left(\left(x_n,y\right),R_n\right)\cap \mathrm{Gr}(f,[a,b])\right)
	&\ge \sum_{i=1}^{\big\lceil\frac{d_{x_n,n}-c_{x_n,n}}{r_{n}}\big\rceil}\mathbb{I}\left(\Delta_{x_n,n,i}\cap L_{f}(y)\right)\cdot\frac{\min\left\lbrace \frac{R_n}{2},\mathrm{Osc}\left(f,\Delta_{x_n,n,i}\right)\right\rbrace }{r_n}\\
	&\ge K_1\cdot \left( \frac{R_{n}}{r_{n}}\right)^{\dim_{\mathrm{A}}L_{f}\left(y\right)}\cdot \frac{\min\left\lbrace \frac{1}{2}R_{n},K_2\cdot R_{n}\right\rbrace }{r_n}\\
	&= C\cdot  \left( \frac{R_{n}}{r_{n}}\right)^{\dim_{\mathrm{A}}L_{f}\left(y\right)+1},
\end{aligned}$$
where $ C=\frac{K_1\cdot\min\left\lbrace 1, 2K_2 \right\rbrace}{2}>0 $.
\end{proof}

Next, let us define
$$ \Lambda:=\left\lbrace (r_{n},R_{n})\ \bigg|\ 0<r_n<R_n<1,\  r_{n}\searrow0,\ R_{n}\searrow0 \ \ \mathrm{and} \ \ \frac{R_{n}}{r_{n}}\nearrow\infty \ \mathrm{as} \ n\to\infty \right\rbrace  $$
and one has $\Lambda_\alpha\subset\Lambda$ for any $\alpha\in(0,1)$. From the definition of the Assouad dimension, Lemma \ref{lem3.4} immediately yields the forthcoming result.

\begin{proposition}\label{pro3.5}
   	Under the precondition of Lemma \ref{lem3.4}, if $ (r_{n},R_{n})\in \Lambda $ for large enough $n$, then
	$$  \dim_{\mathrm{A}}\mathrm{Gr}(f,[a,b])\ge \dim_{\mathrm{A}}L_{f}\left(y\right)+1.  $$
\end{proposition}

\begin{remark}\label{rem3.6}
	 The assumption $\mathrm{(A2)}$ allows a limited number of indices $i \in \left\lbrace 1, 2, \cdots, \big\lceil \frac{d_{x,n} - c_{x,n}}{r_n} \big\rceil \right\rbrace$ to violate the condition (\ref{eq94}) without altering the consequence of the aforementioned proposition. For instance, one may assume that there exist $q_1,q_2\in \mathbb{N}$ such that the condition (\ref{eq94}) holds for all $i=q_1,q_1+1,\cdots, \big\lceil \frac{d_{x,n} - c_{x,n}}{r_n}\big\rceil-q_2$. The reason is that there is a sequence of points $x_n\in L_f(y)$ and a binary sequence $ (r_{n},R_{n})\in \Lambda $ such that
     $$ \begin{aligned}
         \sum_{i=d_1}^{\big\lceil\frac{d_{x_n,n}-c_{x_n,n}}{r_{n}}\big\rceil-d_2}\mathbb{I}\left(\Delta_{x_n,n,i}\cap L_{f}(y)\right)&\ge\#\left\lbrace i\in \mathbb{N} \ \Big|\ \Delta_{x_n,n,i}\cap L_f\left(y\right)\neq\emptyset\right\rbrace-(d_1+d_2)\\
         &\ge K_{1}\cdot \left( \frac{R_n}{r_n}\right)^{\dim_{\mathrm{A}}L_{f}\left(y\right)}-(d_1+d_2)\\
         &\ge \frac{K_{1}}{2}\cdot\left( \frac{R_n}{r_n}\right)^{\dim_{\mathrm{A}}L_{f}\left(y\right)}
     \end{aligned}$$
     for large enough $n$ and then
     $$ \begin{aligned}
N_{r_n}\left(Q\left(\left(x_n,y\right),R_n\right)\cap \mathrm{Gr}(f,[a,b])\right)
	&\ge  \sum_{i=d_1}^{\big\lceil\frac{d_{x_n,n}-c_{x_n,n}}{r_{n}}\big\rceil-d_2}\mathbb{I}\left(\Delta_{x_n,n,i}\cap L_{f}(y)\right)\cdot\frac{\min\left\lbrace \frac{R_n}{2},\mathrm{Osc}\left(f,\Delta_{x_n,n,i}\right)\right\rbrace }{r_n}\\
	&\ge \frac{K_{1}}{2}\cdot\left( \frac{R_{n}}{r_{n}}\right)^{\dim_{\mathrm{A}}L_{f}\left(y\right)}\cdot \frac{\min\left\lbrace \frac{1}{2}R_{n},K_2\cdot R_{n}\right\rbrace }{r_n}\\
	&= C'\cdot  \left( \frac{R_{n}}{r_{n}}\right)^{\dim_{\mathrm{A}}L_{f}\left(y\right)+1},
\end{aligned}$$
where $ C'=\frac{C}{2} $, leading to the result of Proposition \ref{pro3.5}. In fact, the assumptions $\mathrm{(A1)}$ and $\mathrm{(A2)}$ can also be weaken as 
	\begin{equation}\label{eq150}
	    \#\left\lbrace i\in \mathbb{N} \ \Big|\ \Delta_{x,n,i}\cap L_f\left(y\right)\neq\emptyset\ \ \mathrm{and} \ \ \mathrm{Osc}\left(f,\Delta_{x,n,i}\right)\ge K\cdot R_{n}\right\rbrace\ge K\cdot \left(\frac{R_n}{r_{n}}\right)^{\dim_{\mathrm{A}}L_{f}\left(y\right)} 
	\end{equation}  
	for some $ K>0 $.
\end{remark}

In the light of Proposition \ref{pro3.5}, we can now demonstrate Theorem \ref{the1.3} as follows.

\begin{proof}[Proof of Theorem \ref{the1.3}]
    Let
	$$ \rho=\sup_{y\in\mathcal{Y}}\left\lbrace \dim_{\mathrm{A}}L_{f}(y)\right\rbrace. $$
Similar to the proof of Theorem \ref{the1.5}(2), we only need to consider the case when $0<\rho\le1$. Arbitrarily given an $\varepsilon\in(0,\rho)$, there is a $\widetilde{y}\in \mathcal{Y}$ such that
	$$ \dim_{\mathrm{A}}L_f(\widetilde{y})\ge \rho-\varepsilon.$$
    From the definition of $\mathcal{Y}$, there is a binary sequence $(r_{n},R_{n})\in \Lambda_\alpha$
	such that
	$$ \sup_{x\in L_{f}\left(\widetilde{y}\right)}\#\left\lbrace i\in \mathbb{N} \ \Big|\ \Delta_{x,n,i}\cap L_f\left(\widetilde{y}\right)\neq\emptyset\right\rbrace \ge \left( \frac{R_n}{r_{n}}\right)^{\dim_{\mathrm{A}}L_{f}(\widetilde{y})}$$
for all $n\in \mathbb{N}$. Note that $ L_{f}(\widetilde{y}) $ is closed, thus one can find a sequence of points $ \widetilde{x}_{n}\in L_{f}(\widetilde{y}) $ such that
	$$ \begin{aligned}
\#\left\lbrace i\in \mathbb{N} \ \Big|\ \Delta_{\widetilde{x}_{n},n,i}\cap L_f\left(\widetilde{y}\right)\neq\emptyset\right\rbrace
		\ge \frac{1}{2}\cdot\left( \frac{R_{n}}{r_{n}}\right)^{\dim_{\mathrm{A}}L_{f}(\widetilde{y})}
	\end{aligned}   $$
for large enough $n$, which meets the assumption $\mathrm{(A1)}$. Since $f$ is $\alpha$-order count-H\"older continuous and $R_n\le r_n^\alpha$ in the definition of $\Lambda_\alpha$, through a similar argument to (\ref{eq10}) one has
\begin{equation*}
	\begin{aligned}
	\mathrm{Osc}\left(f,\Delta_{\widetilde{x}_{n},n,i}\right)\geq \frac{H^{-}(f)}{2^{\alpha}}\cdot r_n^{\alpha}
	\end{aligned}
\end{equation*}
for all $ i=1,2,\cdots,\big\lceil\frac{d_{x_n,n}-c_{x_n,n}}{r_n}\big\rceil-1$ and large enough $n$, which meets the assumption $\mathrm{(A2)}$ restated in Remark \ref{rem3.6}.
Therefore, Proposition \ref{pro3.5} gives that 
$$ \dim_{\mathrm{A}}\mathrm{Gr}(f,[a,b])\ge \dim_{\mathrm{A}}L_{f}(\widetilde{y})+1\ge \rho+1-\varepsilon. $$
Due to the arbitrariness of $\varepsilon$, we just get the desired result.
\end{proof}

\begin{remark}
   Observe that the precondition of Theorem \ref{the1.3} just serves as a particular case of the assumption $\mathrm{(A2)}$ since the $\alpha$-order count-H\"older continuity of $f$ and $(r_{n},R_{n})\in \Lambda_\alpha$ leads to the formula (\ref{eq94}). Moreover, the aforementioned proof shows that there is a $\widetilde{y}\in \mathcal{Y}$, a sequence of points $ \widetilde{x}_{n}\in L_{f}(\widetilde{y}) $ and a binary sequence $(r_{n},R_{n})\in \Lambda_\alpha$ meeting the condition (\ref{eq150}) where $K$ can be selected as $\min\left\lbrace \frac{1}{4}, \frac{H^{-}(f)}{2^{\alpha}}\right\rbrace$.
\end{remark}

Through an analogous argument to the proof of Theorem \ref{the1.3}, one can also arrive at the following conclusion.

\begin{proposition}\label{cor6.2}
	Let $f\in C_{[a,b]}$ be $\alpha$-order count-H\"older continuous for $\alpha\in(0,1)$. If the following set
	\begin{center}
		$ \begin{aligned}
			\mathcal{Z}:=\Bigg\{ (z_{n})\ \Bigg|\ & z_{n}\in [a,b],\ \mathrm{there}\ \mathrm{is} \ \mathrm{a} \ \xi\in[a,b] \ \mathrm{and}\ \mathrm{a} \ \mathrm{sequence}\ \left( \delta_{n}\right)\\ & \hspace{1cm} \mathrm{such}\ \mathrm{that} \ z_{n}\nearrow \xi,  \ \left(\frac{z_{n+1}-z_{n}}{\delta_{n}}, z_{n+1}-z_{n}\right)\in \Lambda_{\alpha}\
			 \mathrm{and}\\
		& \hspace{2cm}   \#\left\lbrace  i\in \mathbb{N} \ \bigg| \ \widetilde{\Delta}_{n,i}\cap L_{f}\left(f\left(z_{n}\right)\right)\neq \emptyset \right\rbrace\ge \delta_{n}^{\dim_{\mathrm{A}}\Omega_{(z_{n})} }, \ \forall n\in \mathbb{N}   \Bigg\}
		\end{aligned}$
	\end{center}
is non-empty, then 
	$$ \dim_{\mathrm{A}}\mathrm{Gr}(f,[a,b])\ge\sup_{(z_{n})\in\mathcal{Z}}\left\lbrace \dim_{\mathrm{A}}\Omega_{(z_{n})} \right\rbrace+1,   $$
	where $$ \widetilde{\Delta}_{n,i}=\left[z_{n}+ \frac{z_{n+1}-z_{n}}{\delta_{n}}(i-1),z_{n}+\frac{z_{n+1}-z_{n}}{\delta_{n}}i\right]\cap[z_{n},z_{n+1}] $$ 
    for $i=1,2,\cdots,\lceil\delta_n\rceil$ and $$\Omega_{(z_{n})}=\bigcup_{n\in \mathbb{N}} \left\lbrace L_{f}\left(f\left(z_{n}\right)\right)\cap\left[z_{n},z_{n+1}\right]\right\rbrace.    $$
\end{proposition}

The assumption in Proposition \ref{cor6.2} can be viewed as a relaxed version of that in Lemma \ref{lem3.4}, where we just take $R_n=z_{n+1}-z_{n}$ and $r_n=\frac{z_{n+1}-z_{n}}{\delta_{n}}$. Specifically, it no longer requires that the sequence of points $(z_{n})$ must all lie in the same level set. This result reveals that, under certain conditions, the Assouad dimension of the union of countably many level sets—partitioned by some monotonically convergent sequences—plus one, can also provide a lower bound for the Assouad dimension of the graph of a count-H\"older continuous function.

\subsection{Brief explanations for Theorems \ref{the1.7} and \ref{the1.9}} According to the property of the Assouad spectrum introduced earlier in Section \ref{subsection2.2}, Theorem \ref{the1.7} is an immediate consequence from Theorem \ref{the1.5} as long as we take $\theta\to 1^{-}$. It is evident that $$\mathbb{S}_{\mathrm{UB}}\left(\mathrm{Gr}(f,[a,b]),\mathbb{Y}\right)=\sup\limits_{y\in \mathscr{R}_f[a,b]} \left\lbrace \overline{\dim}_{\mathrm{B}}L_f(y)\right\rbrace,$$ 
thus the lower bound of Theorem \ref{the1.9} is gotten. Since the graph of a continuous function on a closed interval is a compact set in $\mathbb{R}^{2}$, we just take $F=\mathrm{Gr}(f,[a,b])$ and $W=\mathbb{Y}$ in Proposition \ref{pro1.1}. Note that $$ \dim_{\mathrm{A}}\mathrm{Proj}_{\mathbb{Y}}\left(\mathrm{Gr}(f,[a,b])\right)=\dim_{\mathrm{A}}\mathscr{R}_{f}[a,b]=1.$$ 
Then by Proposition \ref{pro1.1}, there exists a $ \Omega\in \mathrm{Tan}(\mathrm{Gr}(f,[a,b])) $ and $ x\in \mathrm{Proj}_{\mathbb{Y}}(\Omega) $ such that
$$ {\dim}_{\mathrm{H}}\left( \mathrm{Proj}_{\mathbb{Y}}^{-1}(x)\cap\Omega\right)\ge \dim_{\mathrm{A}}\mathrm{Gr}(f,[a,b])-\dim_{\mathrm{A}}\mathrm{Proj}_{\mathbb{Y}}\left(\mathrm{Gr}(f,[a,b])\right).$$
That is,
$$ \dim_{\mathrm{A}}\mathrm{Gr}(f,[a,b])\le {\dim}_{\mathrm{H}}\left( \mathrm{Proj}_{\mathbb{Y}}^{-1}(x)\cap\Omega\right)+1,$$
which just leads to the upper bound of Theorem \ref{the1.9}. 

Furthermore, it is worth mentioning that Theorem \ref{the1.7} also yields another result
$$ \sup\limits_{y\in \mathscr{R}_f[a,b]} \left\lbrace \overline{\dim}_{\mathrm{B}}L_f(y)\right\rbrace+1-\alpha\le \overline{\dim}_{\mathrm{B}}\mathrm{Gr}(f,[a,b])\le 2-\alpha$$
when letting $\theta\to 0^{+}$. This is exactly consistent with Proposition \ref{rem1.3} that
$$ \overline{\dim}_{\mathrm{B}}\mathrm{Gr}(f,[a,b])=2-\alpha$$
for any $ \alpha $-order H\"older and count-H\"older continuous function $f$ on $[a,b]$ because the upper box dimension of any level set of $f$ cannot exceed 1.

\section{Generalizations to a broader class of H\"older and count-H\"older continuous functions}\label{sec5}

The goal of this section is to extend our previously obtained results to a wider class of functions with the H\"older and count-H\"older continuity, mainly including Weierstrass-type functions, self-affine curves and Brownian motion.

\subsection{Weierstrass-type functions}\label{sec4.1}

In this subsection, we first focus on Weierstrass-type functions $\Phi:\mathbb{R}\to\mathbb{R}$ of the form
\begin{equation}\label{eq5.1}
	\Phi(x):=\sum_{n=0}^{\infty}\beta_{n}^{-\alpha}\phi\left(\beta_{n}x+\theta_{n}\right),
\end{equation} 
which satisfies the following conditions:
\begin{enumerate}
	\item[(\romannumeral1)] $ \alpha\in (0,1) $, $ \beta_{n}>0 $ and $ \theta_{n}\in \mathbb{R} $ for all $ n\in \mathbb{N} $;
	\item[(\romannumeral2)] There are two constants $ \ell_{1}$ and $\ell_{2} $ such that $ 1<\ell_{1}\le\frac{\beta_{n+1}}{\beta_{n}}\le\ell_{2} $ for all $ n\in \mathbb{N} $;
	\item[(\romannumeral3)] $ \phi : \mathbb{R} \to \mathbb{R} $ is a non-constant periodic Lipschitz function.
\end{enumerate}  

\begin{remark}\label{rem5.1}
	In the expression of (\ref{eq5.1}), if we take $\phi(x)=\cos(2\pi x)$, $ \alpha=-\frac{\log\lambda}{\log\beta}$, $\beta_{n}=\beta^{n}$ and $\theta_{n}=0$ for all $n\in \mathbb{N} $ where \( \beta > 1 \) and \( \lambda \in ( \frac{1}{\beta}, 1 ) \), then $\Phi(x)$ 
	turns into the Weierstrass function $W_{\lambda,\beta}(x)$ given in (\ref{eq1}). If we take $\phi(x)=\mathrm{dist}(x,\mathbb{Z})$ and the other parameters are the same, then $\Phi(x)$ becomes the Takagi function $T_{\lambda,\beta}(x)$ given in (\ref{eq2}). In these two cases, we have $\ell_{1}=\ell_{2}=\beta>1$.
\end{remark}

\begin{theorem}\label{the5.1}
	The function $\Phi(x)$ defined as (\ref{eq5.1}) is $ \alpha $-order H\"older and count-H\"older continuous on $ [a,b] $ for large enough $\ell_{1}$.
\end{theorem}

\begin{proof}
	Assume that $ \mathrm{Osc}(\phi,\mathbb{R})=M>0$ and the period of $\phi$ is $T>0$. On one hand, given $ 0<h<\beta_{0}^{-1} $, let $ m $ be the integer satisfying
	\begin{equation}\label{eq6.5}
		\beta^{-1}_{m+1}<h\le\beta^{-1}_{m}. 
	\end{equation}
	Since $\phi$ is Lipschitz on $\mathbb{R}$, there is a certain constant $L>0$ such that
	$$ \begin{aligned}
		\left|\phi\left(\beta_{n}t+\theta_{n}\right)-\phi\left(\beta_{n}u+\theta_{n}\right)\right|
		\le L\beta_{n}\cdot\left|t-u\right|
	\end{aligned}  $$
	for all $ t,u\in \mathbb{R} $ and $ n\in \mathbb{N} $. Combining (\ref{eq6.5}) and $ \frac{\beta_{n+1}}{\beta_{n}}\ge \ell_{1}>1 $ for all $ n\in \mathbb{N} $, it can be estimated that
	$$ \sum_{n=0}^{m}\beta_{n}^{1-\alpha}\le\sum_{j=0}^{m}\ell_{1}^{(\alpha-1)j}\beta_{m}^{1-\alpha}\le\frac{\beta_{m}^{1-\alpha}}{1-\ell_{1}^{\alpha-1}}\le \frac{h^{\alpha-1}}{1-\ell_{1}^{\alpha-1}} $$
	and
	$$ \sum_{n=m+1}^{\infty}\beta_{n}^{-\alpha}\le\sum_{j=0}^{\infty}\ell_{1}^{-\alpha j}\beta_{m+1}^{-\alpha}\le\frac{\beta_{m+1}^{-\alpha}}{1-\ell_{1}^{-\alpha}}\le\frac{h^{\alpha}}{1-\ell_{1}^{-\alpha}}. $$
	Consequently, we obtain
	\begin{equation}\label{eq12}
		\begin{aligned}
			\left|\Phi(x+h)-\Phi(x)\right|
			&\le\sum_{n=0}^{m}\beta_{n}^{-\alpha}\left|\phi\left(\beta_{n}(x+h)+\theta_{n}\right)-\phi\left(\beta_{n}x+\theta_{n}\right)\right|\\
			&\quad+\sum_{n=m+1}^{\infty}\beta_{n}^{-\alpha}\left|\phi\left(\beta_{n}(x+h)+\theta_{n}\right)-\phi\left(\beta_{n}x+\theta_{n}\right)\right|\\
			&\le Lh\cdot\sum_{n=0}^{m} \beta_{n}^{1-\alpha}+M\cdot\sum_{n=m+1}^{\infty}\beta_{n}^{-\alpha}\\
			&\le\frac{Lh\cdot h^{\alpha-1}}{1-\ell_{1}^{\alpha-1}}+\frac{M\cdot  h^{\alpha}}{1-\ell_{1}^{-\alpha}}\\
			&\le H^{+}(\Phi)\cdot h^{\alpha},
		\end{aligned} 
	\end{equation}
	where
	$$ H^{+}(\Phi)=\frac{L}{1-\ell_{1}^{\alpha-1}}+\frac{M}{1-\ell_{1}^{-\alpha}}>0. $$
	This implies that $ \Phi(x) $ is an $ \alpha $-order H\"older continuous function on $ [a,b] $. 
	
	On the other hand, for all $ m\in \mathbb{N} $ and $ 0<h\le T\beta^{-1}_{m} $, one can show that
	\begin{equation}\label{eq3.4}
		\begin{aligned}
			&\ \beta_{m}^{-\alpha}\left|\phi\left(\beta_{m}(x+h)+\theta_{m}\right)-\phi\left(\beta_{m}x+\theta_{m}\right)\right|-\left|\Phi(x+h)-\Phi(x)\right|\\ 
			&\le\left| \Phi(x+h)-\Phi(x)-\beta_{m}^{-\alpha}\left[ \phi\left(\beta_{m}(x+h)+\theta_{m}\right)-\phi\left(\beta_{m}x+\theta_{m}\right)\right] \right| \\
			&=\sum_{n=0}^{m-1}\beta_{n}^{-\alpha}\left|\phi\left(\beta_{n}(x+h)+\theta_{n}\right)-\phi\left(\beta_{n}x+\theta_{n}\right)\right| \\
			&\quad+\sum_{n=m+1}^{\infty}\beta_{n}^{-\alpha}\left|\phi\left(\beta_{n}(x+h)+\theta_{n}\right)-\phi\left(\beta_{n}x+\theta_{n}\right)\right|\\
			&\le\dfrac{L\ell_{1}^{\alpha-1}T\cdot \beta_{m}^{-\alpha}}{1-\ell_{1}^{\alpha-1}}+\frac{M\ell_{1}^{-\alpha}\cdot\beta_{m}^{-\alpha}}{1-\ell_{1}^{-\alpha}}\\
			&=\left( \dfrac{L\cdot T}{\ell_{1}^{1-\alpha}-1}+\dfrac{M}{\ell_{1}^{\alpha}-1}\right)\cdot \beta_{m}^{-\alpha}
		\end{aligned} 
	\end{equation}
	through a similar discussion to (\ref{eq12}). Write $ \delta_{0}=T\beta_{1}^{-1}  $. For each $ 0<\delta\le\delta_{0} $, let $ m_{0} $ be the integer such that
	\begin{equation}\label{eq3.5}
		T\beta_{m_{0}}^{-1}<\delta\le T\beta_{m_{0}-1}^{-1}.  
	\end{equation}
	For each $ x\in[a,b] $, there exists an $ h_{0}>0 $ such that $ \beta_{m_{0}}h_{0}\in (0,T] $ and then
	\begin{equation}\label{eq3.6}
		\left|\phi\left(w+\beta_{m_{0}}h_{0}\right)-\phi\left(w\right)\right|\ge \frac{M}{2},
	\end{equation}
	where $ w=\beta_{m_{0}}x+\theta_{m_{0}} $. At this time, one has
    $$ 0<h_{0}\le T\beta_{m_{0}}^{-1}<\delta. $$
    Therefore, together with the assumption $ \frac{\beta_{n+1}}{\beta_{n}}\le \ell_{2} $ for all $ n\in \mathbb{N} $, it holds by (\ref{eq3.4})--(\ref{eq3.6}) that
	$$ \begin{aligned}
		\left|\Phi\left(x+h_{0}\right)-\Phi(x)\right|&\ge 		\beta_{m_{0}}^{-\alpha}\left|\phi\left(\beta_{m_{0}}(x+h_{0})+\theta_{m_{0}}\right)-\phi\left(\beta_{m_{0}}x+\theta_{m_{0}}\right)\right| -\left( \dfrac{L\cdot T}{\ell_{1}^{1-\alpha}-1}+\dfrac{M}{\ell_{1}^{\alpha}-1}\right)\cdot\beta_{m_{0}}^{-\alpha}\\
		&\ge \ell_{2}^{-\alpha}\left( \frac{M}{2}-\dfrac{L\cdot T}{\ell_{1}^{1-\alpha}-1}-\dfrac{M}{\ell_{1}^{\alpha}-1}\right)\cdot\beta_{m_{0}-1}^{-\alpha}\\
		&\ge H^{-}(\Phi)\cdot\delta^{\alpha},
	\end{aligned} $$
	where
	\begin{equation}\label{eq78}
	    H^{-}(\Phi)=\frac{1}{\left(T\ell_{2}\right)^{\alpha}}\cdot\left( \frac{M}{2}-\dfrac{L\cdot T}{\ell_{1}^{1-\alpha}-1}-\dfrac{M}{\ell_{1}^{\alpha}-1}\right)>0
	\end{equation}  
	for large enough $ \ell_{1} $. Hence, $ \Phi(x) $ is an $ \alpha $-order count-H\"older continuous function on $ [a,b] $, completing the proof.
\end{proof}

\begin{remark}
	It can be seen from the above proof that the precondition of $ \frac{\beta_{n+1}}{\beta_{n}}\le \ell_{2} $ for all $n\in \mathbb{N}$ is indispensable for determining the count-H\"older continuity of $ \Phi(x) $. In other words, if $ \frac{\beta_{n+1}}{\beta_{n}}\nearrow\infty $ as $n\to\infty$, one cannot judge the count-H\"older continuity of such functions like the Besicovitch function with rapidly growing frequencies (see \cite{Baranski2012}, for instance) whose graphs may have the different lower and upper box dimensions.
\end{remark}

On the basis of Remark \ref{rem5.1}, Proposition \ref{pro2.6} can exactly be considered as a direct consequence of Theorem \ref{the5.1} in fact. Furthermore, Theorems \ref{the1.5} and \ref{the1.9} immediately yield the undermentioned generalized result with Proposition \ref{cor1.3} turning into a particular case.      

\begin{proposition}
	Provided $\ell_{1}$ is large enough, it holds that
	\begin{equation}\label{eq5.2}
		\dim_\mathrm{A}^{\theta}\mathrm{Gr}(\Phi,[a,b])\leq \begin{cases}
			1+\displaystyle\frac{1-\alpha}{1-\theta}, \quad & \theta\in(0,\alpha),\\ \\
			2,\quad & \theta\in[\alpha,1),
		\end{cases}  
	\end{equation} 
	
	\begin{equation}\label{eq5.2New}
		\dim_\mathrm{A}^{\theta}\mathrm{Gr}(\Phi,[a,b])\ge \begin{cases}
			\sup\limits_{y\in \mathscr{R}_\Phi[a,b]} \left\lbrace \overline{\dim}_{\mathrm{B}}L_\Phi(y)\right\rbrace+\dfrac{1-\alpha}{1-\theta},\quad & \theta\in(0,\alpha),\\
			\\
			\sup\limits_{y\in \mathscr{R}_\Phi[a,b]} \left\lbrace \overline{\dim}_{\mathrm{B}}L_\Phi(y)\right\rbrace+1,\quad &\theta\in[\alpha,1)
		\end{cases}
	\end{equation}  
	and    
	$$ \begin{aligned}
		\mathbb{S}_{\mathrm{UB}}\left(\mathrm{Gr}(\Phi,[a,b]),\mathbb{Y}\right)+1\le\dim_{\mathrm{A}}\mathrm{Gr}(\Phi,[a,b])
		\le\sup\limits_{\Omega\in\mathrm{Tan}(\mathrm{Gr}(\Phi,[a,b]))} \mathbb{S}_{\mathrm{H}}\left(\Omega,\mathbb{Y}\right)+1,
	\end{aligned}$$
	where $\Phi(x)$ has been defined as (\ref{eq5.1}).
\end{proposition}

From the above discussion, one can provide a variety of H\"older and count-H\"older continuous functions with the same Assouad spectrum range as (\ref{eq5.2}) and (\ref{eq5.2New}). Unlike the classical Weierstrass and Takagi functions, the ratio $\frac{\beta_{n+1}}{\beta_n}$ can be selected to meet the condition
$$ \ell_{1}<\liminf_{n\to\infty}\frac{\beta_{n+1}}{\beta_{n}}<\limsup_{n\to\infty}\frac{\beta_{n+1}}{\beta_{n}}<\ell_{2}$$
rather than being restricted to a constant value. In addition, the function $ \phi(x) $ can be taken as any other non-constant periodic Lipschitz functions like $ \sin x $, $ \cos^{2}x $ or the ‘zig-zag’ function of period $ \gamma>0 $ given by
$$ Z(\gamma k+t)=\begin{cases}
	t,\quad &t\in \left[0,\dfrac{\gamma}{4}\right),\\\\
	\dfrac{\gamma}{2}-t,\quad & t\in \left[\dfrac{\gamma}{4},\dfrac{3\gamma}{4} \right),\\\\
	t-\gamma,\quad & t\in \left[\dfrac{3\gamma}{4},\gamma\right),
\end{cases} $$
where $ k\in\mathbb{N}$ and $ t\in[0,\gamma) $.

Another typical and noteworthy example is the following Weierstrass function which is appended a random phase to each term
\begin{equation*}
	W_{\lambda,\beta,\Theta}(x) = \sum_{n=0}^{\infty} \lambda^n \cos(2\pi (\beta^n x+\theta_{n})), 
\end{equation*}
where $ \Theta=\left\lbrace \theta_{0},\theta_{1},\cdots,\theta_{n},\cdots \right\rbrace  $ and $\theta_{n}$ is a random number within $[0,1]$. The graph of $ W_{\lambda,\beta,\Theta}$ can be interpreted as a fractal basin boundary for a random dynamical system \cite{Yu2024} and Hunt \cite{Hunt1998,Hunt1991} has shown that its Hausdorff dimension is equal to the box dimension $ 2+\frac{\log\lambda}{\log\beta}$ if each $\theta_{n}$ is selected independently with respect to the uniform probability measure in $[0,1]$. Now we also derive its Assouad dimension by reason of its H\"older and count-H\"older continuity if $\beta$ is large enough.

\subsection{Self-affine functions}

Next, this subsection turns towards certain self-affine sets formed by iterated function systems(IFS), which can also be graphs of functions as long as we choose suitable affine transformations. Let us consider a system of affine transformations \(\{\mathscr{W}_i\}_{i=1}^N\) where	each \(\mathscr{W}_i : \mathbb{R}^2 \to \mathbb{R}^2\) has the form
	\[
	\mathscr{W}_i \begin{pmatrix} x \\ y \end{pmatrix} = 
	\begin{pmatrix} \alpha_i & 0 \\ a_i & c_i \end{pmatrix} 
	\begin{pmatrix} x \\ y \end{pmatrix} + 
	\begin{pmatrix} \beta_i \\ b_i \end{pmatrix},
	\]
where \(\alpha_i \in (0,1)\) is the horizontal contraction ratio, \(c_i \in (-1,1) \setminus \{0\}\) is the vertical contraction ratio, \(\beta_i = \sum_{j=1}^{i-1} \alpha_j\) ensures domain coverage, and \(a_i, b_i \in \mathbb{R}\) are translation parameters. Now we define 
$$ 		\mathcal{P}_1= \mathscr{W}_1(	\mathcal{P}_1) = \left(0, \frac{b_1}{1 - c_1}\right)\quad \mathrm{and} \quad 
			\mathcal{P}_N= \mathscr{W}_N(	\mathcal{P}_N) = \left(1, \frac{a_N + b_N}{1 - c_N}\right) $$
as the fixed points of $ \mathscr{W}_1 $ and $ \mathscr{W}_N $. 
	Consider an IFS $ \mathscr{I}=\left\lbrace [0,1]\times\mathbb{R}; 	\mathscr{W}_i,i=1,2,\cdots,N\right\rbrace $ with the following conditions:
	\begin{enumerate}
		\item[(\romannumeral1)] (Domain partition) \(\sum_{i=1}^N \alpha_i = 1\) with \(\alpha_i >0 \) for all $ i=1,2,\cdots,N $;
			\item[(\romannumeral2)] (Vertical dominance) $ \alpha_i \le \left| c_{i}\right| $  for all $ i=1,2,\cdots,N $;
		\item[(\romannumeral3)] (Continuity condition) \(\mathscr{W}_i(\mathcal{P}_N) = \mathscr{W}_{i+1}(\mathcal{P}_1)\) for all \(i=1,2,\cdots,N\);
		\item[(\romannumeral4)] (Non-degeneracy) The points \(\{\mathscr{W}_1(\mathcal{P}_1), \mathscr{W}_2(\mathcal{P}_1), \cdots, \mathscr{W}_N(\mathcal{P}_1), \mathcal{P}_N\}\) are not collinear.
	\end{enumerate}
	Then, the IFS $ \mathscr{I}$ admits a unique attractor which is the graph of a continuous function \(\varphi: [0,1] \to \mathbb{R}\) satisfying
 \begin{equation}\label{eq87}
 	\mathrm{Gr}(\varphi,[0,1])=\bigcup_{i=1}^{N}\mathscr{W}_i\left(\mathrm{Gr}(\varphi,[0,1])\right).  
 \end{equation}
	The upcoming theorem indicates the H\"older and count-H\"older continuity of $ \varphi $.
	\begin{theorem}\label{the4.2}
		The function $ \varphi(x) $ defined above is $ \frac{\log c_{\max}}{\log \alpha_{\min}} $-order H\"older continuous and $ \frac{\log c_{\min}}{\log \alpha_{\max}} $-order count-H\"older continuous on $ [0,1] $ where
		$$ \alpha_{\max}=\max_{1\le i\le N}\alpha_{i},  \quad \alpha_{\min}=\min_{1\le i\le N}\alpha_{i},  $$
		$$  c_{\max}=\max_{1\le i\le N} |c_i|, \quad c_{\min}=\min_{1\le i\le N} |c_i|. $$ 
	\end{theorem}

\begin{proof}
Let \(\Sigma_k = \{1,2,\cdots,N\}^k\) be the set of multi-indices of length \(k\in\mathbb{N}_+\). For each $$ \mathbf{i}_{k} := (i_1,i_2,\cdots,i_k) \in \Sigma_k, $$ we define the composite map \(\mathscr{W}_{\mathbf{i}_{k}} := \mathscr{W}_{i_1} \circ \mathscr{W}_{i_2}\circ\cdots \circ \mathscr{W}_{i_k}\), the vertical contraction factor \(C_{\mathbf{i}_{k}}: =\prod_{j=1}^{k}c_{i_j} \), and the horizontal projection \(I_{\mathbf{i}_{k}}: = \mathscr{W}_{\mathbf{i}_{k}}([0,1])\) with the length
$ \alpha_{\mathbf{i}_{k}}:= \prod_{j=1}^{k}\alpha_{i_j}  $. The condition (\romannumeral2) implies that $ \alpha_{\mathbf{i}_{k}}\le \left|C_{\mathbf{i}_{k}}\right|  $.
	Based on the above assumptions, the linear part of \(\mathscr{W}_{\mathbf{i}_{k}}\) can be represented by the following matrix
	\[
	\mathscr{V}_{\mathbf{i}_{k}} := \begin{pmatrix} \alpha_{\mathbf{i}_{k}} & 0 \\ A_{\mathbf{i}_{k}} & C_{\mathbf{i}_{k}} \end{pmatrix},
	\]
	where 
	\[
	A_{\mathbf{i}_{k}} = \sum_{j=1}^k a_{i_j} \left( \prod_{\ell=1}^{j-1} c_{i_\ell} \right) \left( \prod_{\ell=j+1}^{k} \alpha_{i_\ell} \right).
	\]
Since $ 0<c_{\max}<1 $, it can be estimated that
	$$ \left|A_{\mathbf{i}_{k}}\right|\le a_{\max} \sum_{j=1}^kc_{\max} ^{j-1}\frac{\alpha_{\mathbf{i}_{k}}}{\alpha_{i_{j}}}\le \frac{a_{\max}}{\alpha_{\min}}\sum_{j=0}^{k-1}c_{\max} ^{j}\alpha_{\mathbf{i}_{k}}\le K_{1}\alpha_{\mathbf{i}_{k}}\le K_{1}\left| C_{\mathbf{i}_{k}}\right|,  $$
	where $ a_{\max} = \max_{1\le i\le N} |a_i| $ and $K_{1} = \frac{a_{\max}}{\alpha_{\min} (1 - c_{\max})}$. Observe that $ \mathrm{Gr}(\varphi,I_{\mathbf{i}_{k}}) $ is the
	affine image $ \mathscr{W}_{\mathbf{i}_{k}}\left(\mathrm{Gr}(\varphi,[0,1])\right)  $ which can be considered as a translate of $ 	\mathscr{V}_{\mathbf{i}_{k}}\left( \mathrm{Gr}(\varphi,[0,1])\right)  $. 
	Let $ H $ be the height of $ \mathrm{Gr}(\varphi,[0,1]) $ and then the image $ 	\mathscr{V}_{\mathbf{i}_{k}}\left( \mathrm{Gr}(\varphi,[0,1])\right)  $ is contained in a rectangle of height $ (K_{1}+H)\left| C_{\mathbf{i}_{k}}\right| $. Additionally, it is trivial from the condition (\romannumeral4) of non-degeneracy that \(\mathscr{W}_{\mathbf{i}_{k}}(\mathrm{Gr}(\varphi,[0,1]))\) contains a triangle with height no less than \( K_{2} \left|C_{\mathbf{i}_{k}}\right|\) for some \(K_{2} > 0\). Therefore, the range of \(\varphi\) over \(I_{\mathbf{i}_{k}}\) satisfies 
	\begin{equation}\label{eq:range}
		K_{2} \left|C_{\mathbf{i}_{k}}\right| \leq \mathrm{Osc}(\varphi,I_{\mathbf{i}_{k}}) \leq (K_{1}+H) \left|C_{\mathbf{i}_{k}}\right|.
	\end{equation}

On one hand, given $ 0<\delta\le 1 $, there is an integer $ k_{1} $ such that
$$  \alpha_{\min}^{k_{1}+1}\le\delta< \alpha_{\min}^{k_{1}}.  $$
Because $ \frac{\delta}{ \alpha_{\min}}\ge\alpha_{\min}^{k_{1}} $, one has
$$ k_{1}\log c_{\max}\le  \frac{\log c_{\max}}{\log \alpha_{\min}}\cdot\log\frac{\delta}{ \alpha_{\min}}, $$
which yields that
\begin{equation}\label{eq95}
	 c_{\max}^{k_{1}}\le \left(\frac{\delta}{ \alpha_{\min}}\right)^{ \frac{\log c_{\max}}{\log \alpha_{\min}}}. 
\end{equation}
Combined with 
 $$  \left|I_{\mathbf{i}_{k_{1}}}\right|=\prod_{j=1}^{k_{1}}\alpha_{i_j}\ge \alpha_{\min}^{k_{1}}>\delta,\quad \forall  \mathbf{i}_{k_{1}}\in \Sigma_{k_{1}},   $$
for each $ x\in [0,1] $ there are two adjacent intervals $ I_{\mathbf{i}'_{k_{1}}},I_{\mathbf{i}''_{k_{1}}}\in \Sigma_{k_{1}} $ such that $[x-\delta,x+\delta]\subset I_{\mathbf{i}'_{k_{1}}}\cup I_{\mathbf{i}''_{k_{1}}} $. Thus, it follows from (\ref{eq:range}) and (\ref{eq95}) that 
$$\begin{aligned}
	\mathrm{Osc}\left(\varphi,[x-\delta,x+\delta]\cap[0,1]\right)&\le\mathrm{Osc}\left(\varphi,I_{\mathbf{i}'_{k_{1}}}\cup I_{\mathbf{i}''_{k_{1}}}\right)\\ 
	&\le (K_{1}+H) \left( \left|C_{\mathbf{i}'_{k_{1}}}\right|+\left|C_{\mathbf{i}''_{k_{1}}}\right|\right)\\
	&\le  2(K_{1}+H)c_{\max}^{k_{1}} \\
	&\le H^{+}(\varphi)\cdot\delta^{\frac{\log c_{\max}}{\log \alpha_{\min}}},
\end{aligned} $$
where $ H^{+}(\varphi)= 2(K_{1}+H)\alpha_{\min}^{-\frac{\log c_{\max}}{\log \alpha_{\min}}}>0 $. This leads to the H\"older continuity of $ \varphi $ on $ [0,1] $ with the order $ \frac{\log c_{\max}}{\log \alpha_{\min}} $.

On the other hand, given $ 0<\delta\le 1 $, there is an integer $ k_{2} $ such that
$$  \alpha_{\max}^{k_{2}}<\delta\le \alpha_{\max}^{k_{2}-1}.  $$
Note that $\delta\alpha_{\max}\le \alpha_{\max}^{k_{2}} $, so
$$  \frac{\log c_{\min}}{\log \alpha_{\max}}\cdot\log\left( \delta\alpha_{\max}\right) \le  k_{2}\log c_{\min},  $$
which gives
 \begin{equation}\label{eq96}
 	 c_{\min}^{k_{2}}\ge \left( \delta\alpha_{\max}\right)^{\frac{\log c_{\min}}{\log \alpha_{\max}}}.
 \end{equation}
Since $$  \left|I_{\mathbf{i}_{k_{2}}}\right|=\prod_{j=1}^{k_{2}}\alpha_{i_j}\le \alpha_{\max}^{k_{2}}<\delta,\quad \forall  \mathbf{i}_{k_{2}}\in \Sigma_{k_{2}},   $$
for each $ x\in [0,1] $ there is an $ \mathbf{i}'_{k_{2}}\in \Sigma_{k_{2}} $ such that $ I_{\mathbf{i}'_{k_{2}}}\subset[x-\delta,x+\delta]\cap[0,1] $. 
Consequently, by (\ref{eq:range}) and (\ref{eq96}), 
$$\mathrm{Osc}\left(\varphi,[x-\delta,x+\delta]\cap[0,1]\right)\ge\mathrm{Osc}\left(\varphi,I_{\mathbf{i}'_{k_{2}}}\right) \ge 	K_{2} \left|C_{\mathbf{i}_{k_{2}}}\right|\ge K_{2} c_{\min}^{k_{2}}\ge H^{-}(\varphi)\cdot\delta^{\frac{\log c_{\min}}{\log \alpha_{\max}}}, $$
where $ H^{-}(\varphi)=K_{2}\alpha_{\max}^{\frac{\log c_{\min}}{\log \alpha_{\max}}}>0 $. This just means that $ \varphi $ is a count-H\"older continuous function on $ [0,1] $ with the order $ \frac{\log c_{\min}}{\log \alpha_{\max}} $.
\end{proof}

In the light of Theorems \ref{the1.5} and \ref{the1.9}, one can come to the following conclusion for the Assouad dimension and spectrum of $ \mathrm{Gr}(\varphi,[0,1]) $.

\begin{proposition}
	For the function $\varphi(x)$ defined as (\ref{eq87}), it holds that
		\begin{equation*}
		\dim_\mathrm{A}^{\theta}\mathrm{Gr}(\varphi,[0,1])\leq \begin{cases}
			1+\displaystyle\frac{1-\frac{\log c_{\max}}{\log \alpha_{\min}}}{1-\theta}, \quad & \theta\in\left(0,\dfrac{\log c_{\max}}{\log \alpha_{\min}}\right),\\ \\
			2,\quad & \theta\in\left[\dfrac{\log c_{\max}}{\log \alpha_{\min}},1\right),
		\end{cases}  
	\end{equation*} 
	
	\begin{equation*}
		\dim_\mathrm{A}^{\theta}\mathrm{Gr}(\varphi,[0,1])\ge \begin{cases}
			\sup\limits_{y\in \mathscr{R}_\varphi[0,1]} \left\lbrace \overline{\dim}_{\mathrm{B}}L_\varphi(y)\right\rbrace+\dfrac{1- \frac{\log c_{\min}}{\log \alpha_{\max}}}{1-\theta},\quad & \theta\in\left(0, \dfrac{\log c_{\min}}{\log \alpha_{\max}}\right),\\
			\\
			\sup\limits_{y\in \mathscr{R}_\varphi[0,1]} \left\lbrace \overline{\dim}_{\mathrm{B}}L_\varphi(y)\right\rbrace+1,\quad &\theta\in\left[\dfrac{\log c_{\min}}{\log \alpha_{\max}},1\right)
		\end{cases}
	\end{equation*}  
	and    
	$$ \begin{aligned}
		\mathbb{S}_{\mathrm{UB}}\left(\mathrm{Gr}(\varphi,[0,1]),\mathbb{Y}\right)+1\le\dim_{\mathrm{A}}\mathrm{Gr}(\varphi,[0,1])
		\le\sup\limits_{\Omega\in\mathrm{Tan}(\mathrm{Gr}(\varphi,[0,1]))} \mathbb{S}_{\mathrm{H}}\left(\Omega,\mathbb{Y}\right)+1,
	\end{aligned}$$
where $ \alpha_{\max} $, $ \alpha_{\min} $, $ c_{\max} $ and $ c_{\min} $ have been given earlier in Theorem \ref{the4.2}.
\end{proposition}

\begin{remark}
As a particular case of uniform horizontal contractions, if we choose \(\alpha_i = \frac{1}{N}\) for all \(i=1,2,\cdots,N\), then $ \varphi(x) $ defined as (\ref{eq87}) is H\"older continuous with the order $ -\frac{\log c_{\max}}{\log N} $ and count-H\"older continuous with the order $ -\frac{\log c_{\min}}{\log N} $ on $ [0,1] $. Ulteriorly, the Assoaud spectrum of $ \mathrm{Gr}(\varphi,[0,1]) $ will possess the following estimation:
	\begin{equation*}
	\dim_\mathrm{A}^{\theta}\mathrm{Gr}(\varphi,[0,1])\leq \begin{cases}
		1+\displaystyle\frac{1+\frac{\log c_{\max}}{\log N}}{1-\theta}, \quad & \theta\in\left(0,-\dfrac{\log c_{\max}}{\log N}\right),\\ \\
		2,\quad & \theta\in\left[-\dfrac{\log c_{\max}}{\log N},1\right)
	\end{cases}  
\end{equation*} 
and
\begin{equation*}
	\dim_\mathrm{A}^{\theta}\mathrm{Gr}(\varphi,[0,1])\ge \begin{cases}
		\sup\limits_{y\in \mathscr{R}_\varphi[0,1]} \left\lbrace \overline{\dim}_{\mathrm{B}}L_\varphi(y)\right\rbrace+\dfrac{1+\frac{\log c_{\min}}{\log N}}{1-\theta},\quad & \theta\in\left(0, -\dfrac{\log c_{\min}}{\log N}\right),\\
		\\
		\sup\limits_{y\in \mathscr{R}_\varphi[0,1]} \left\lbrace \overline{\dim}_{\mathrm{B}}L_\varphi(y)\right\rbrace+1,\quad &\theta\in\left[-\dfrac{\log c_{\min}}{\log N},1\right).
	\end{cases}
\end{equation*} 
\end{remark}

\begin{remark}[Fractal interpolation functions]
    In fact, those foregoing results for self-affine functions are also applicable to fractal interpolation functions which can be constructed through arbitrary points \(\{(x_i, y_i)\}_{i=0}^N\) with non-uniform spacing. More precisely, set
\[
\alpha_i = x_i - x_{i-1}, \quad \beta_i = x_{i-1}
\]
and choose suitable \(a_i, b_i, c_i\) to satisfy the continuity at given interpolation points while controlling the H\"older and count-H\"older continuity via the vertical contraction ratios. On the basis of it, the corresponding Assouad dimension and spectrum can be derived for graphs of fractal interpolation functions as well.
\end{remark}

\subsection{Brownian motion}\label{subsec4.3}

At the end of this section, we add some remarks on Brownian motion (or the Wiener process) $ \mathcal{B}(t) $, which serves as a common example of L\'evy processes \cite{Levy1934} and is defined to be
the stochastic process with the following properties:
	\begin{enumerate}
		\item[(1)] \(\mathcal{B}(0) = 0\);
		\item[(2)] (Independence of increments) For all \(0< t_1 < t_2 < \cdots < t_n\), the random variables \(\mathcal{B}(t_2) - \mathcal{B}(t_1)\), \(\mathcal{B}(t_3) - \mathcal{B}(t_2)\), \(\cdots\), \(\mathcal{B}(t_n) - \mathcal{B}(t_{n-1})\) are independent;
		\item[(3)] (Gaussian increments) For all $ t,h>0 $, \(\mathcal{B}(t+h) - \mathcal{B}(t)\) is
		normally distributed with mean 0 and variance $ h $;
		\item[(4)] (Continuity) \(\mathcal{B}(t)\) is almost surely a continuous function of $ t $.
	\end{enumerate}

Geometric properties, such as various fractal dimensions, of the graph of Brownian motion have been a
well investigated field for a long time. While its Hausdorff dimension was first calculated by Taylor \cite{Taylor1953} to be $ \frac{3}{2} $ almost surely, Howroyd and Yu \cite{Howroyd2019} have further claimed that it almost surely possesses the full Assouad dimension of 2, which can be a remarkable progress. Apart from its dimensional results, our most concern for the analytical properties of Brownian motion is that it almost surely has the $ \lambda $-order H\"older continuity on $ [0,1] $ for each $ \lambda\in(0,\frac{1}{2}) $ (see \cite[Proposition 16.1]{Falconer2014}), which yields from the first claim of Theorem \ref{the1.5} that almost surely
	\begin{equation*}
	\dim_\mathrm{A}^{\theta}\mathrm{Gr}(\mathcal{B},[0,1])\leq \begin{cases}
		1+\displaystyle\frac{1-\lambda}{1-\theta}, \quad & \theta\in\left(0,\lambda\right),\\ \\
		2,\quad & \theta\in\left[\lambda,1\right).
	\end{cases}  
\end{equation*} 

From the aforementioned analysis, one can notice that obtaining a lower bound for the Assouad dimension hinges on establishing a uniform lower bound for the oscillation of the corresponding function over every sufficiently small interval—a condition guaranteed by any count-H\"older continuous function. Although the count-H\"older continuity of Brownian motion cannot be determined globally on the entire domain, a weaker result regarding its amplitude over infinitely many intervals has been established as follows:

\begin{proposition}[\cite{Howroyd2019}]\label{pro4.3}
For $ \frac{1}{2}\le \lambda\le 1 $, almost surely there is a sequence of closed intervals $ I_{n}=[a_n,b_n]\subset [0,1] $ and a sequence of integers $ N_{n}\nearrow \infty $ as $ n\to\infty $ such that
	$$ \mathrm{Osc}\left(\mathcal{B}, \left[a_{n}+(i-1)\frac{\left|I_{n}\right|}{N_{n}},a_{n}+i\frac{\left|I_{n}\right|}{N_{n}}\right]\right)\ge \left|I_{n}\right|^{\lambda}  $$
	for all $ i=1,2,\cdots,N_{n} $ and all $ n\in \mathbb{N} $.
\end{proposition}

Here we omit the detailed proof of Proposition \ref{pro4.3} because it can be regarded as a direct consequence from the argument in \cite[Section 3]{Howroyd2019} where the authors derived it by the $ \frac{1}{2} $-scaling property and independent increments of Brownian motion together with the Borel-Cantelli lemma. Applying Proposition \ref{pro4.3} and Lemma \ref{lem2.1}, in the next section, we will provide an alternative explanation for why Brownian motion has the full Assouad dimension of 2, as also demonstrated in \cite{Howroyd2019}.
	
	\section{Discussions on the full Assouad dimension of graphs}\label{sec6}

  As a particularly meaningful case, the problem of characterizing continuous functions with graphs of the full Assouad dimension has received significant attention. This section will deal with certain circumstances under which the Assouad dimension of graphs can reach 2.

\subsection{Consequences of main results}
	
In the first place, we give some pivotal consequences in this subsection. The upcoming proposition follows directly from Theorem \ref{the1.7}.

\begin{proposition}\label{cor3}
	Let $f$ be an $ \alpha $-order count-H\"older continuous function on $[a,b]$ for $\alpha\in (0,1)$. If there is a level set of $f$ whose upper box dimension is 1, then $\dim_{\mathrm{A}}\mathrm{Gr}(f,[a,b])=2$.
\end{proposition}

Proposition \ref{cor3} explicitly demonstrates that identifying a level set with the full upper box dimension is a valid strategy for establishing that the graph of a count-Hölder continuous function also has the full Assouad dimension. Conversely, if no level set attains an upper box dimension of 1, the next proposition following from Theorem \ref{the1.9} offers an alternative approach for determining whether, under certain conditions, some level set may achieve the full Assouad dimension.

\begin{proposition}\label{cor4}
     Let $f$ be an $ \alpha $-order count-H\"older continuous function on $[a,b]$ for $\alpha\in (0,1)$. If there is a $y\in\mathscr{R}_{f}[a,b]$ and a binary sequence $(r_{n},R_{n})\in \Lambda_{\alpha}$ such that
	$$	\sup_{x\in L_{f}(y)} \#\left\lbrace  i\in \mathbb{N} \ \Big| \ \Delta_{x,n,i}\cap L_{f}(y)\neq \emptyset \right\rbrace\ge \frac{R_{n}}{r_{n}}$$
    for large enough $n$, then $\dim_{\mathrm{A}}\mathrm{Gr}(f,[a,b])=2$.
\end{proposition}

In reality, the function $f\in C_{[a,b]}$ does not necessarily possess the count-H\"older continuity on $[a,b]$. In the following, we state a general result which is derived from Proposition \ref{pro3.5}.

\begin{proposition}\label{pro4}
      Let $f\in C_{[a,b]}$. Suppose that there is a $ y\in \mathscr{R}_{f}[a,b] $, a sequence of points $x_n\in L_f(y)$ and a binary sequence $(r_{n},R_{n})\in \Lambda_{\alpha}$ satisfying both assumptions $\mathrm{(A1)}$ and $\mathrm{(A2)}$ for large enough $n$. If $\dim_{\mathrm{A}}L_f(y)=1$, then $\dim_{\mathrm{A}}\mathrm{Gr}(f,[a,b])=2$.
\end{proposition}

Based on this premise, it results directly from assumption $\mathrm{(A1)}$ that $\dim_{\mathrm{A}}L_{f}(y)=1$.
We emphasize that a comparison of Propositions \ref{cor4} and \ref{pro4} leads to the following observation: To achieve the objective, if $f$ is $\alpha$-order count-H\"older continuous, the required binary sequence $(r_{n}, R_{n})$ must belong to the set $\Lambda_{\alpha}$. However, if $f$ satisfies the stronger assumption $\mathrm{(A2)}$, then $(r_{n}, R_{n})$ only needs to meet a weaker criterion, that is, $(r_{n}, R_{n})\in\Lambda$ is sufficient.

\subsection{A class of constructive examples with the full Assouad dimension}

Building on the preceding propositions, this subsection will provide a category of explicit examples of continuous functions whose graphs can attain the full Assouad dimension. 

First, we present some sequences with the full Assouad dimension. To facilitate the subsequent discussion, given a $ \zeta\in \mathbb{R} $, let $ \Pi(\zeta) $ denote the family of real number sequences $ \left( u_{n}\right)_{\mathbb{N}_+}$ which converge to the only cluster point $ \zeta $ and satisfy the following properties:
\begin{enumerate}
	 \item[(\romannumeral1)] $\left|u_{n+1}-u_{n}\right|\searrow 0^+ $ as $n\to\infty$;
     \item[(\romannumeral2)] $\frac{u_{n+1}-\zeta}{u_{n}-\zeta}\nearrow 1^-$ as $n\to\infty$.
\end{enumerate}
We set $\widehat{\Omega}=\left\lbrace u_{n}  \ \big|\  
 n\in\mathbb{N}_+\right\rbrace$. Note that the property (\romannumeral2) implies either $\zeta<u_{n+1}<u_{n}$ or $u_{n}<u_{n+1}<\zeta$ for all $n\in \mathbb{N}_+$. Thus, $u_{n}$ must converge monotonically to $\zeta$ as $n\to\infty$. Additionally, since both the sequence in $ \Pi(\zeta) $ and its step length exhibit monotonic behavior, we refer to such sequence as bi-monotonic. 

\begin{example}
     Some examples of usually encountered sequences belonging to $ \Pi(0) $ have been given below:\\ (1) $\frac{1}{\log\log(n+2)}$; (2) $\frac{1}{\log(n+1)}$; (3) $\frac{1}{n^{\alpha}}$ $(\alpha>0)$;  (4) $\frac{1}{n^{\alpha}\log^{\beta}(n+1)}$ $(\alpha,\beta>0)$; (5) $n^{-\frac{1}{\sqrt{\log (n+1)}}}$; (6) $n^{-\frac{1}{\log \log (n+1)}}$.\\
     Here $n\in \mathbb{N}_+$.
\end{example}

\begin{remark}
  We point out that geometric sequences $ a^n$ with $0<a<1$ are excluded from $ \Pi(0) $ owing to their failure to satisfy the property (\romannumeral2).
\end{remark}

Fixing $k\ge1$ and $ \left( u_{n}\right)_{\mathbb{N}_+}\in  \Pi(\zeta) $, we now define $ \widehat{r}_{n}=k\left| u_{n+1}-u_{n}\right| $, $\widehat{R}_{n}=2\left|u_{n}-\zeta\right|$, $ \widehat{c}_{n}=\min\left\lbrace u_{n},\zeta \right\rbrace  $, $ \widehat{d}_{n}=\max\left\lbrace u_{n},\zeta \right\rbrace$ and $$ \widehat{\Delta}_{n,i}=\left[ \widehat{c}_{n}+ (i-1) \widehat{r}_{n}, \widehat{c}_{n}+i \widehat{r}_{n}\right]\cap\left[ \widehat{c}_{n}, \widehat{d}_{n}\right] $$ for $ i=1,2,\cdots,\Big\lceil\frac{\widehat{d}_{n}-\widehat{c}_{n}}{\widehat{r}_{n}}\Big\rceil$. For these definitions of notations, the undermentioned result is fundamental.

\begin{proposition}\label{pro56}
   Provided $n$ is large enough, we have the following assertions. 
   \begin{enumerate}
       \item[(1)]  $\left(\widehat{r}_{n},\widehat{R}_{n}\right)\in \Lambda$;
       \item[(2)] $  \#\left\lbrace i\in \mathbb{N} \ \Big|\  \widehat{\Delta}_{n,i}\cap \widehat{\Omega}\neq\emptyset\right\rbrace\ge\dfrac{1}{2}\cdot\dfrac{\widehat{R}_{n}}{\widehat{r}_{n}}     $.
   \end{enumerate}
\end{proposition}

\begin{proof}
    (1) In the light of the properties (\romannumeral1) and (\romannumeral2) of $ \left( u_{n}\right)_{\mathbb{N}_{+}}\in \Pi(\zeta) $, one has $\widehat{r}_{n}\searrow 0$ and $ \widehat{R}_{n}\searrow 0  $ as $n\to\infty$. Furthermore, it can be deduced that
   $$ \frac{\widehat{R}_{n}}{\widehat{r}_{n}}=\frac{2}{k}\left|\frac{u_{n}-\zeta}{ u_{n+1}-u_{n}}\right|=\frac{\frac{2}{k}}{1-\frac{u_{n+1}-\zeta}{u_{n}-\zeta}}\nearrow\infty,\quad n\to\infty.$$
   Consequently, it gives
   $ 0<\widehat{r}_{n}<\widehat{R}_{n}<1$ for large enough $n$, which leads to the required result.

   (2)  The property (\romannumeral1) of $ \left( u_{n}\right)_{\mathbb{N}_{+}}\in \Pi(\zeta) $ ensures that $\left| u_{p+1}-u_{p}\right|\le\widehat{r}_{n}  $ for all $ p\in \mathbb{N} $ and $ p\ge n $. Hence, for each $ i=1,2,\cdots,\Big\lceil\frac{\widehat{d}_{n}-\widehat{c}_{n}}{\widehat{r}_{n}}\Big\rceil $, there exists at least a $ p_{i}\in \mathbb{N} $ and $ p_{i}\ge n $ such that $ u_{p_{i}}\in\widehat{\Delta}_{n,i}  $, which just means that $ \widehat{\Delta}_{n,i}\cap \widehat{\Omega}\neq\emptyset  $. Observe that $\widehat{d}_{n}-\widehat{c}_{n}=\left|u_{n}-\zeta\right|=\frac{\widehat{R}_{n}}{2}$. Therefore,
   $$ \#\left\lbrace i\in \mathbb{N} \ \Big|\  \widehat{\Delta}_{n,i}\cap \widehat{\Omega}\neq\emptyset\right\rbrace=\Bigg\lceil\frac{\widehat{d}_{n}-\widehat{c}_{n}}{\widehat{r}_{n}}\Bigg\rceil \ge\frac{1}{2}\cdot\dfrac{\widehat{R}_{n}}{\widehat{r}_{n}}, $$
   finishing the present proof.
\end{proof}

The foregoing proposition also implies the set $\widehat{\Omega}$ has the full Assouad dimension of 1, which has been stated previously in the first claim of Proposition \ref{pro1.3}. Below we give a brief explanation of it for completeness.

\begin{proof}[Proof of Proposition \ref{pro1.3}(1)]
 With Proposition \ref{pro56}, there are two sequences of scales $0<\widehat{r}_{n}<\widehat{R}_{n}<1$ for large enough $n$ such that
  $$   N_{\widehat{r}_{n} }\left(Q\left(\zeta,\widehat{R}_{n}\right)\cap\widehat{\Omega}\right)\ge N_{\widehat{r}_{n} }\left(\left[ \widehat{c}_{n}, \widehat{d}_{n}\right]\cap\widehat{\Omega}\right)=\#\left\lbrace i\in \mathbb{N} \ \Big|\  \widehat{\Delta}_{n,i}\cap \widehat{\Omega}\neq\emptyset\right\rbrace\ge\dfrac{1}{2}\cdot\dfrac{\widehat{R}_{n}}{\widehat{r}_{n}}\nearrow\infty,\quad n\to\infty. $$
  Then, it yields that $\dim_{\mathrm{A}}\widehat{\Omega}=1$ from the definition of the Assouad dimension.
\end{proof}

\begin{proposition}\label{pro53}
		Fix $k\ge 1$, $ \left( u_{n}\right)_{\mathbb{N}_{+}}\in \Pi(\zeta) $. Let $ f\in C_{[a,b]} $ satisfy $f(u_{n})=f(\zeta)$ for all $n\in\mathbb{N}_{+}$. If there is a constant $ L>0  $ and a $q\in \mathbb{N}$ such that 
		\begin{equation}\label{eq75}
		    \mathrm{Osc}\left(f,\widehat{\Delta}_{n,i}\right)\ge L\cdot \left|u_{n}-\zeta\right| 
		\end{equation} 
		for all $ i=q,q+1,\cdots,\Big\lceil\left|\frac{u_{n}-\zeta}{k (u_{n+1}-u_{n})}\right|\Big\rceil - 2 $ and for large enough $ n$, then   $\dim_{\mathrm{A}}\mathrm{Gr}(f,[a,b])=2  $.
\end{proposition}

\begin{proof}
    	Under this premise, it is trivial that $ \widehat{\Omega}\cup \left\lbrace \zeta\right\rbrace\subseteq L_f(\zeta) $ and then
	$$ \dim_{\mathrm{A}}\widehat{\Omega}\le \dim_{\mathrm{A}}L_f(\zeta)\le 1.  $$
  Since $\dim_{\mathrm{A}}\widehat{\Omega}=1$, one has $\dim_{\mathrm{A}}L_f(\zeta)=1$. For large enough $n$, Proposition \ref{pro56} allows that there exists an $f(\zeta)\in \mathscr{R}_{f}[a,b]$, a sequence of points $\widehat{c}_{n}\in L_f(f(\zeta))$ and a binary sequence $\left(\widehat{r}_{n},\widehat{R}_{n}\right)\in \Lambda$ such that
     $$ \#\left\lbrace i\in \mathbb{N} \ \Big|\  \Delta_{\widehat{c}_{n},n,i}\cap L_f(f(\zeta))\neq\emptyset\right\rbrace\ge\#\left\lbrace i\in \mathbb{N} \ \Big|\  \widehat{\Delta}_{n,i}\cap \widehat{\Omega}\neq\emptyset\right\rbrace\ge\dfrac{1}{2}\cdot\left(\dfrac{\widehat{R}_{n}}{\widehat{r}_{n}}\right)^{\dim_{\mathrm{A}}L_f(\zeta)}, $$
    which meets the assumption $\mathrm{(A1)}$. Besides, the condition (\ref{eq75}) ensures that there exist certain $q_1,q_2\in \mathbb{N}$ such that
    $$  \mathrm{Osc}\left(f,\Delta_{\widehat{c}_{n},n,i}\right)\ge L\cdot \left|u_{n}-\zeta\right|=\frac{L}{2}\cdot \widehat{R}_{n}$$
  for all $i=q_1,q_1+1,\cdots, \Big\lceil \frac{d_{\widehat{c}_{n},n} - c_{\widehat{c}_{n},n}}{\widehat{r}_{n}}\Big\rceil-q_2$, which meets the assumption $\mathrm{(A2)}$  restated in Remark \ref{rem3.6}. Together with Proposition \ref{pro4}, we just get our desired result.
\end{proof}

With the above theoretical preparations, we will construct a class of continuous functions in the subsequent parts.
Given a decreasing sequence $ \left( u_{n}\right)_{\mathbb{N}_{+}}\in \Pi(\zeta) $, let us define the function $ \Psi:[\zeta,u_1]\to \mathbb{R} $ by
\begin{equation}\label{eq3.5}
	 \Psi(x)=\begin{cases}
	     \psi(x)\cdot(\mathcal{F}\circ\mathcal{G})(x), & x\in (\zeta,u_1],\\ \\
         0, & x=\zeta,
	 \end{cases}
\end{equation}
where the functions $\psi$, $\mathcal{F}$ and $\mathcal{G}$ are defined as follows:
\begin{enumerate}
 \item[(1)] $\psi : [\zeta,u_1]\to \mathbb{R} $ is a continuous and strictly increasing function satisfying $\psi(\zeta)=0$ and $\psi(u_n-u_{n+1})\ge u_n-\zeta $ for large enough $n$;
    \item[(2)] $\mathcal{G}:(\zeta,u_1]\to \mathbb{R}$ is a continuous and strictly decreasing function satisfying $ \lim_{x\to \zeta^{+}}\mathcal{G}(x)=+\infty$;
    \item[(3)] $\mathcal{F}:\mathbb{R}\to \mathbb{R}$ is a non-constant periodic Lipschitz function and $\left\lbrace \mathcal{G}(u_{n}) \ \big| \ n\in \mathbb{N}_{+}\right\rbrace  $ is the set of all zero points of $ \mathcal{F} $ on $ [\mathcal{G}(u_1),+\infty) $;
\end{enumerate}
The function $\Psi$ constructed as above is a continuous function on $ [\zeta,u_1] $ with a single unbounded variation point $x=\zeta$, which satisfies $\widehat{\Omega}\cup\left\lbrace \zeta \right\rbrace=L_{\Psi}(0)$. Namely, $\Psi(u_n)=\Psi(\zeta)=0$ for all $n\in \mathbb{N}_+$. Additionally, one can come to the following conclusion that its graph on $[\zeta,u_1]$ has the full Assouad dimension.

\begin{proposition}\label{pro5.5}
  $\dim_{\mathrm{A}}\mathrm{Gr}(\Psi,[\zeta,u_1])=2$.
\end{proposition}

\begin{proof}
   Assume that $\max_{x\in \mathbb{R}}\left|\mathcal{F}(x)\right|=M>0$ and the period of $\mathcal{F}$ is $T>0$. Let $k_0$ be a positive integer such that
   \begin{equation}\label{eq9}
       \mathcal{G}(u_{n+k_0})-\mathcal{G}(u_{n})=T
   \end{equation} 
   for all $ n\in \mathbb{N}_{+} $. At this time, we set
   $$\widetilde{\Delta}_{n,i}=\left[\zeta+(k_0+1)(u_{n}-u_{n+1})(i-1), \zeta+(k_0+1)(u_{n}-u_{n+1})i\right]\cap[\zeta,u_n]$$ 
   for all $ i=1,2,\cdots,\big\lceil\frac{u_{n}-\zeta}{(k_0+1)(u_{n}-u_{n+1})}\big\rceil $. Then for all $ p\in \mathbb{N}_{+} $ and $ p\ge n $, 
   $$ \frac{(k_0+1)(u_{n}-u_{n+1})}{u_{p}-u_{p+k_0+1}}> \frac{(k_0+1)(u_{n}-u_{n+1})}{\sum_{j=n+1}^{n+k_0+1}u_{j-1}-u_{j}}>1,$$
   which means that there are at least $k_0+1$ zero points of $\mathcal{F}\circ\mathcal{G}$ within each $\widetilde{\Delta}_{n,i}$ for $ i=1,2\cdots,\big\lceil\frac{u_{n}-\zeta}{(k_0+1)(u_{n}-u_{n+1})}\big\rceil -1 $. Thus by (\ref{eq9}), there is an $x_{n,i}\in \widetilde{\Delta}_{n,i}$ such that $\left|(\mathcal{F}\circ\mathcal{G})(x_{n,i})\right|=M$. Since $\psi $ is strictly increasing on $[\zeta,u_1]$, for all $x\in \widetilde{\Delta}_{n,i}$, we can claim that 
   $$ \psi(x)\ge\psi((k_0+1)(u_{n}-u_{n+1}))>\psi(u_{n}-u_{n+1})\ge u_n-\zeta $$
   for all $ i=2,3,\cdots,\big\lceil\frac{u_{n}-\zeta}{(k_0+1)(u_{n}-u_{n+1})}\big\rceil -2 $ for large enough $n$.
   Then, it can be estimated that
    $$\mathrm{Osc}\left(\Psi,\widetilde{\Delta}_{n,i}\right)\ge \left|\Psi(x_{n,i})\right|=\psi(x_{n,i})\cdot \left|(\mathcal{F}\circ\mathcal{G})(x_{n,i})\right|\ge M\cdot (u_n-\zeta),$$
    which exactly meets the condition (\ref{eq75}). Combined with $ \Psi\in C_{[\zeta,u_1]} $ and $\Psi(u_n)=\Psi(\zeta)=0$ for all $n\in \mathbb{N}_+$, this result can be a direct consequence of Proposition \ref{pro53}. 
\end{proof}

For any $f\in C_{[a,b]}$ and $t,h\in \mathbb{R}$, it is simple to verify that
\begin{equation}\label{eq91}
    \dim_{\mathrm{A}}\mathrm{Gr}(g,[t-b,t-a])=\dim_{\mathrm{A}}\mathrm{Gr}(f,[a,b]),
\end{equation} 
where $g(x):=f(-x+t)+h$ is obtained by applying some translation and symmetry transformations to the function $f$.
 From the construction of (\ref{eq3.5}) with Proposition \ref{pro5.5}, the second claim of Proposition \ref{pro1.3} holds trivially. Here, we provide a concise explication of it.

\begin{proof}[Proof of Proposition \ref{pro1.3}(2)]
Fix $h\in \mathbb{R}$ and $ \left( u_{n}\right)_{\mathbb{N}_{+}}\in \Pi(\zeta) $. 
If $u_n\searrow\zeta$ as $n\to\infty$, we define the function $\Psi_1:[\zeta,u_1]\to\mathbb{R}$ by
$$\Psi_1(x):=\Psi(x)+h. $$
If $u_n\nearrow\zeta$ as $n\to\infty$, we define the function $\Psi_2:[u_1,\zeta]\to\mathbb{R}$ by
$$ \Psi_2(x):=\Psi(-x+\zeta+u_1)+h.$$
By the property of $\Psi$, one has $\Psi_1\in C_{[\zeta,u_1]}$, $\Psi_2\in C_{[u_1,\zeta]}$ and $\Psi_1(u_n)=\Psi_1(\zeta)=\Psi_2(u_n)=\Psi_2(\zeta)=h$ for all $n\in \mathbb{N}_+$. Meanwhile, combining Proposition \ref{pro5.5} and (\ref{eq91}),
$$ \dim_{\mathrm{A}}\mathrm{Gr}(\Psi_1,[\zeta,u_1])=\dim_{\mathrm{A}}\mathrm{Gr}(\Psi_2,[u_1,\zeta])=\dim_{\mathrm{A}}\mathrm{Gr}(\Psi,[\zeta,u_1])=2,$$
which shows the existence of our required continuous functions.
\end{proof}

At the end of this subsection, we present a specific example of the function $\Psi$ when $\zeta=0$ and $u_1=1$.

\begin{example}\label{ex5.9}
Let
    $$\widetilde{\Psi}(x)=\begin{cases}
        x^{\alpha}\cdot\sin\dfrac{\pi}{x^{\beta}},& x\in(0,1],\\
        \\
        0, & x=0,
    \end{cases}$$
    where $\beta>0$ and $0<\alpha<\frac{1}{1+\beta}$. Then $\dim_{\mathrm{A}}\mathrm{Gr}\left(\widetilde{\Psi},[0,1]\right)=2$.
\end{example}

\begin{remark}
    In the above example, we take $u_n=n^{-\frac{1}{\beta}}$, $\psi(x)=x^{\alpha}$, $\mathcal{G}(x)=\frac{1}{x^{\beta}}$ and $\mathcal{F}(x)=\sin(\pi x)$ for $\alpha,\beta>0$. If $\alpha<\frac{1}{1+\beta}$, it can be estimated that $$\left(n^{-\frac{1}{\beta}}-(n+1)^{-\frac{1}{\beta}}\right)^\alpha\ge n^{-\frac{1}{\beta}} $$ for large enough $n$. Thus, $(u_n)_{\mathbb{N}_+}$, $\psi$, $\mathcal{G}$ and $\mathcal{F}$ satisfy all the given assumptions, which ensures that $\mathrm{Gr}\left(\widetilde{\Psi},[0,1]\right)$ has the full Assouad dimension.
\end{remark}

\begin{figure}[htbp] 
	\centering 
	\begin{subfigure}[b]{0.48\textwidth} 
		\centering
		\includegraphics[width=3.2in, height=2.2in]{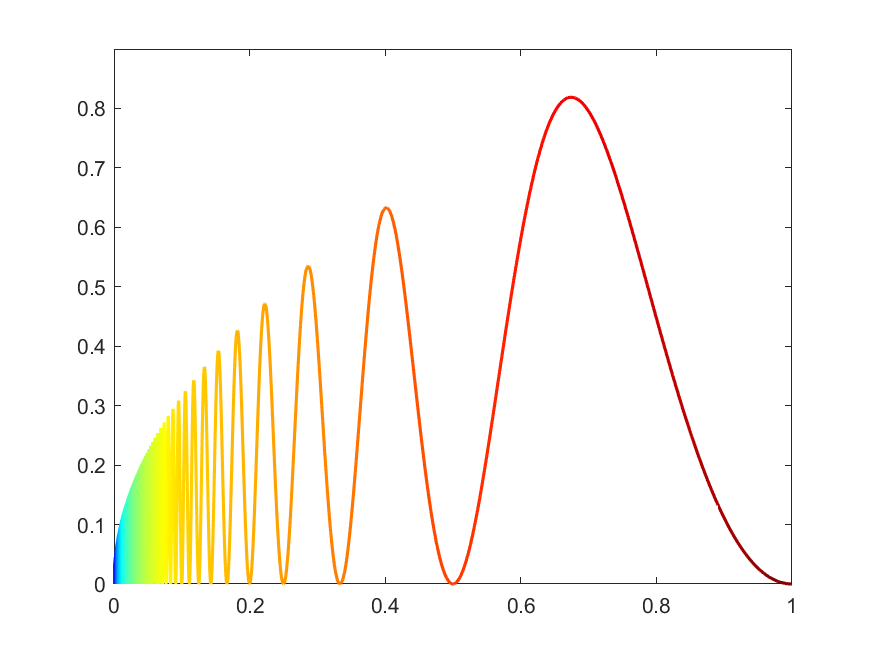} 
		\caption{The graph of $y=\sqrt{x}\cdot\sin^2\dfrac{\pi}{x}$ on $[0,1]$.} 
		\label{fig:sub1}
	\end{subfigure}
	\hfill 
	\begin{subfigure}[b]{0.48\textwidth} 
		\centering
		\includegraphics[width=3.2in, height=2.2in]{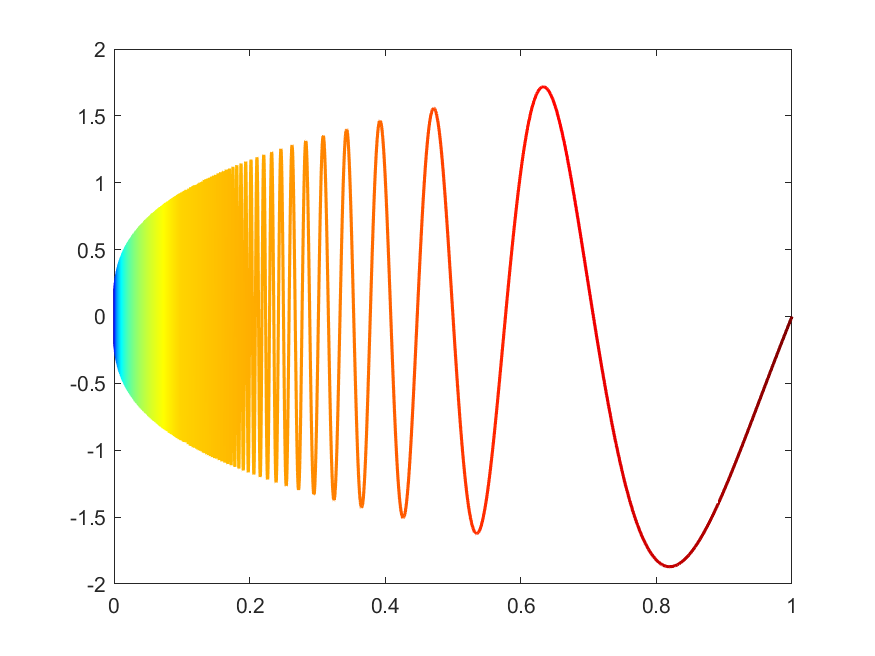} 
		\caption{The graph of $y=2x^{\frac{1}{3}}\cdot\sin\dfrac{\pi}{x^2}$ on $[0,1]$.} 
		\label{fig:sub2}
	\end{subfigure}
	\caption{Graphs of two examples of $\widetilde{\Psi}$ with the full Assouad dimension.}\label{fig5} 
	\label{fig:parallel}
\end{figure}

The crucial point in giving the concrete form of (\ref{eq3.5}) lies in the determination of $\psi$ and $u_n$, requiring that they should satisfy $\psi(u_n-u_{n+1})\ge u_n-\zeta $. In fact, if $\alpha=\frac{1}{1+\beta}$ in Example \ref{ex5.9}, we can choose $\psi(x)=\gamma\cdot x^{\alpha}$ for $\gamma\ge \beta^\alpha$, which also meets this condition. Furthermore, $\mathcal{F}$ can be taken as any other non-constant periodic Lipschitz functions like $ \cos x $, $ \sin^{2}x $ or the ‘zig-zag’ function discussed earlier in Section \ref{sec4.1}, and $\mathcal{G}$ is accordingly adjusted to be a suitable function. Figure \ref{fig5} visually shows the graphs of two special examples of $\widetilde{\Psi}$.

\subsection{Applications of Lemma \ref{lem2.0}} Throughout this paper, Lemma \ref{lem2.0} serves as the primary technical tool for estimating the lower bound of the Assouad dimension of graphs of continuous functions. We obtain the following lemma by applying it directly.

\begin{lemma}\label{lem52}
	Let $ f\in C_{[a,b]} $. Suppose that there is a constant $ K>0 $, a sequence of points $t_n\in [a,b]$ and two sequences of scales $0<r_{n}<R_{n}<1$ such that the following two assumptions hold for large enough $ n $:
	\begin{enumerate}
		\item[$ (\mathrm{H1}) $]  $   \#\left\lbrace i\in \mathbb{N} \ \big|\ \Delta_{n,i}\cap L_f\left(f(t_{n})+\frac{2\lambda-1}{2}R_{n}\right)\neq\emptyset\right\rbrace\ge K\cdot \frac{d_{n}-c_{n}}{r_{n}}   $ for some $ \lambda\in(0,1) $;
		\item[$ (\mathrm{H2}) $] $ \mathrm{Osc}\left(f,\Delta_{n,i}\right)\ge K\cdot R_{n} $ for all $ i=1,2,\cdots,\lceil\frac{d_{n}-c_{n}}{r_{n}}\rceil$.
	\end{enumerate}
	Then
		$$ 	 \begin{aligned}
		N_{r_{n}}\left(Q\left((t_{n},f(t_{n})),R_{n}\right)\cap \mathrm{Gr}(f,[a,b])\right)
		\ge C'\cdot\left(\frac{R_{n}}{r_{n}}\right)^{2},
	\end{aligned}
	$$
	where $ C' $ is a certain positive constant, $ c_{n}=\max\left\lbrace t_{n}-\frac{R_{n}}{2},a \right\rbrace  $, $ d_{n}=\min\left\lbrace t_{n}+\frac{R_{n}}{2},b \right\rbrace  $ and $ \Delta_{n,i}=[c_{n} + (i-1)r_{n},c_{n}+ir_{n}]\cap[c_{n},d_{n}] $ for $ i=1,2,\cdots,\lceil\frac{d_{n}-c_{n}}{r_{n}}\rceil$.
\end{lemma}

\begin{proof}
	Note that $ d_{n}-c_{n}\ge \frac{R_{n}}{2}  $ for large enough $ n$. Then, the assumption $ (\mathrm{H1}) $ ensures that there exists a $ \lambda\in(0,1) $ such that 
	$$ \begin{aligned}
		\sum_{i =1}^{\lceil\frac{d_{n}-c_{n}}{r_{n}}\rceil}\mathbb{I}\left(\Delta_{n,i}\cap L_f\left(f(t_{n})+\frac{2\lambda-1}{2}R_{n}\right)\right)
		=\#\left\lbrace i\in \mathbb{N} \ \bigg|\ \Delta_{n,i}\cap L_f\left(f(t_{n})+\frac{2\lambda-1}{2}R_{n}\right)\neq\emptyset\right\rbrace
		\ge \frac{K}{2}\cdot \frac{R_{n}}{r_{n}}.
	\end{aligned}   $$
Employing Lemma \ref{lem2.0}, combined with the assumption $ (\mathrm{H2}) $, we have
		$$ 	 \begin{aligned}
		&	N_{r_{n}}\left(Q\left((t_{n},f(t_{n})),R_{n}\right)\cap \mathrm{Gr}(f,[a,b])\right)\\
		&\ge \sum_{i =1}^{\lceil\frac{d_{n}-c_{n}}{r_{n}}\rceil}\mathbb{I}\left(\Delta_{n,i}\cap L_f\left(f(t_{n})+\frac{2\lambda-1}{2}R_{n}\right)\right)\cdot \frac{\min\left\lbrace \frac{1-\left|2\lambda-1\right| }{2} R_{n}, \mathrm{Osc}\left(f,\Delta_{n,i}\right)\right\rbrace }{r_{n}}\\
		&\ge \frac{K}{2}\cdot \frac{R_{n}}{r_{n}}\cdot \frac{\min\left\lbrace \frac{1-\left|2\lambda-1\right| }{2} R_{n},K\cdot R_{n}\right\rbrace }{r_{n}}\\
		&\ge C'\cdot \left(\frac{R_{n}}{r_{n}}\right)^{2},
	\end{aligned}
	$$
	where $ C'=\frac{K}{4}\min\left\lbrace \left( 1-\left|2\lambda-1\right|\right) ,2K \right\rbrace  $, completing the present proof.
\end{proof}

\begin{remark}\label{rem5.5}
	Similar to Remark \ref{rem3.6}, the assumptions $ (\mathrm{H1}) $ and $ (\mathrm{H2}) $ in the above lemma can be weaken as 
	$$ \#\left\lbrace i\in \mathbb{N} \ \bigg|\ \Delta_{n,i}\cap L_f\left(f(t_{n})+\frac{2\lambda-1}{2}R_{n}\right)\neq\emptyset\ \ \mathrm{and} \ \ \mathrm{Osc}\left(f,\Delta_{n,i}\right)\ge K\cdot R_{n}\right\rbrace\ge K\cdot \frac{d_{n}-c_{n}}{r_{n}}  $$
	for some $ K>0 $ and $ \lambda\in(0,1) $ without altering the consequence.
\end{remark}

By the definition of the Assouad dimension, Lemma \ref{lem52} immediately yields the forthcoming result.

\begin{proposition}\label{pro5.6}
		Under the precondition of Lemma \ref{lem52}, if $ (r_{n},R_{n})\in \Lambda $ for large enough $ n$, we have
	$  \dim_{\mathrm{A}}\mathrm{Gr}(f,[a,b])=2  $.
\end{proposition}

Compared to Lemma \ref{lem3.4}, the assumption in Lemma \ref{lem52} permits the parameter $\lambda$ to be an arbitrarily given number within $(0,1)$. Together with Remark \ref{rem5.5}, if we choose $(r_{n},R_{n})$ as certain particular binary sequences, one can also acquire a consequence of Proposition \ref{pro5.6} as follows.

\begin{corollary}\label{cor5.14}
    Let $f\in C_{[a,b]}$. If there is a constant $K>0$, a sequence $(\delta_{n})$ and a monotonically increasing convergent sequence of points $(z_{n})$ satisfying $z_n\in [a,b]$, $\left(\frac{z_{n+1}-z_{n}}{\delta_{n}}, z_{n+1}-z_{n}\right)\in \Lambda$ and
    $$ \#\left\lbrace  i\in \mathbb{N} \ \bigg| \ \widetilde{\Delta}_{n,i}\cap L_{f}\left(f\left(z_{n}\right)+(2\lambda-1)(z_{n+1}-z_{n})\right)\neq \emptyset \ \ \mathrm{and} \ \ \mathrm{Osc}\left(f,\widetilde{\Delta}_{n,i}\right)\ge K\cdot (z_{n+1}-z_{n}) \right\rbrace\ge K\cdot\delta_{n}$$
    for infinitely many $n\in \mathbb{N}$ and some $\lambda\in(0,1)$, then 
	$ \dim_{\mathrm{A}}\mathrm{Gr}(f,[a,b])=2 $.
	Here $$ \widetilde{\Delta}_{n,i}=\left[z_{n}+ \frac{z_{n+1}-z_{n}}{\delta_{n}}(i-1),z_{n}+\frac{z_{n+1}-z_{n}}{\delta_{n}}i\right]\cap[z_{n},z_{n+1}] $$ 
    for $i=1,2,\cdots,\lceil\delta_n\rceil$.
\end{corollary}

\begin{remark}
    This corollary represents a special case of Proposition \ref{pro5.6} where we can choose $t_n=z_n$, $R_n=2(z_{n+1}-z_{n})$ and $r_n=\frac{2(z_{n+1}-z_{n})}{\delta_{n}}$. At this time, $\frac{d_n-c_n}{r_n}=\frac{R_n}{r_n}=\delta_n$ for large enough $n$. Moreover, when fixing $\lambda=\frac{1}{2}$, it can also be a corollary following from Proposition \ref{cor6.2} under the condition that $\dim_{\mathrm{A}}\Omega_{(z_{n})}=1$.
\end{remark}

\begin{remark}
    As a representative example, Corollary \ref{cor5.14} just offers a more specific explanation of why Brownian motion $ \mathcal{B}(t) $ introduced in Section \ref{subsec4.3} attains the full Assouad dimension based on \cite{Howroyd2019}. Actually, it can be inferred from the argument of \cite[Theorem 3.1]{Howroyd2019} that there are infinitely many intervals $I_n=[a_n,a_{n+1}]\subset[0,1]$ (we can choose $a_1=0$ and $a_{n+1}=a_n+\frac{1}{2^n}$ for all $n\ge1$) and a sequence of integer $N_n\nearrow\infty$ such that
 $$   \left[a_{n}+(i-1)\frac{b_n-a_n}{N_{n}},a_{n}+i\frac{b_n-a_n}{N_{n}}\right]\cap L_{\mathcal{B}}\left(\mathcal{B}\left(a_{n}\right)+\frac{a_{n+1}-a_{n}}{4}\right)\neq \emptyset$$
 and
 $$ \mathrm{Osc}\left(\mathcal{B}, \left[a_{n}+(i-1)\frac{b_n-a_n}{N_{n}},a_{n}+i\frac{b_n-a_n}{N_{n}}\right]\right)\ge b_n-a_n $$
for all $i=1,2,\cdots,N_{n}$. 
In this case, we just take $z_n=a_n$, $\delta_n=2N_n$ and $\lambda=\frac{3}{4}$, which satisfy all the required conditions stated in Corollary \ref{cor5.14}. Thus by Corollary \ref{cor5.14}, one has $\dim_{\mathrm{A}}\mathrm{Gr}(\mathcal{B},[0,1])=2$.
\end{remark}

\section{Remarks and further research}\label{sec7}

In this final section, we add some remarks on this work and outline future research directions.

\subsection{Sharpness of the lower bound}
We emphasize that the lower bound in Theorem \ref{the1.7} is sharp with respect to level sets, where the upper box dimension cannot be extended to the Assouad dimension. The proof of Theorem 1.1, which establishes the lower bound of the Assouad spectrum for $\alpha$-order count-H\"{o}lder continuous functions, shows that for a level set with the upper box dimension \( s \), there exists a sequence \( r_n\searrow 0 \) as $n\to\infty$ such that inequality (\ref{eq19}) holds. Then, the lower bound for the Assouad dimension can be obtained via the limit behavior of the Assouad spectrum. However, for a level set with the Assouad dimension \( s \), one cannot guarantee the existence of a binary sequence \( (R_n, r_n) \in \Lambda \) that satisfies \( R_n < r_n^\alpha \) such that the assumption (A1) holds for some sequence of points within the level set. This requirement is necessary because, for a count-H\"{o}lder continuous function \( f \), we can only ensure that the oscillation of \( f \) on each interval \( \Delta_i \) of sufficiently small length \( r_n \) is at least of the order \( r_n^\alpha \), but not necessarily of the order \( R_n \). If this weaker oscillation condition fails, then by the covering number estimate in Lemma \ref{lem2.0}, the term \( \min\left\{R_n, \text{Osc}(f, \Delta_i)\right\} \) cannot be guaranteed to be bounded below by a quantity of order \( R_n \). As a result, the total number of required covering boxes will not be bounded below by \( ( \frac{R_n}{r_n} )^{s+1} \).

\subsection{Boundedness of the count-H\"older continuity}

The derivation of our main results relies crucially on the count-H\"older continuity of functions under consideration. Through our investigation in Section \ref{sec5}, these results can be applicable to a wider class of Weierstrass-type functions $\Phi$ of the form (\ref{eq5.1}), which notably include two classical functions of particular interest: the Weierstrass and Takagi functions. However, the proof of Theorem \ref{the5.1} clearly indicates that the count-H\"older continuity of $\Phi$ can only be established when $\ell_1>1$ is large enough for a given $\alpha$ so that the count-H\"older coefficient $H^{-}(\Phi)$ computed via (\ref{eq78}) can become a positive constant. As two special cases, the Weierstrass and Takagi functions can be determined to have the count-H\"older continuity only when $\beta>\frac{1}{\lambda}$ is large enough for a given $\lambda<1$. Therefore, it remains unverified whether the lower bound in Theorem \ref{the1.7} continues to hold for the Weierstrass and Takagi functions when $\beta$ is small—particularly values close to 1.

\subsection{Open questions}

Based on the foregoing discussions, several open questions that could be further considered in the future have been gathered as follows:

\begin{question}
    If $\beta$ is sufficiently close to 1, do the Weierstrass function defined in (\ref{eq1}) and Takagi functions defined in (\ref{eq2}) still exhibit the count-Hölder continuity? In other words, is the lower bound established in Theorem \ref{the1.7} still valid for them?
\end{question}

\begin{question}
    Is there any count-H\"older continuous function having a level set with the full upper box dimension of 1?  More generally, can we find more examples of fractal functions whose graphs have the full Assouad dimension of 2?
\end{question}

\begin{question}
  Can we derive a similar upper bound for the Assouad dimension of graphs of H\"older and count-H\"older continuous functions in terms of level sets? Further, can we obtain an exact value of it?
\end{question}

\section*{Acknowledgments}
The authors thank Jonathan M. Fraser for his valuable discussions and comments, and Balázs Bárány and Roope Anttila for identifying an error in Lemma 3.2 and suggesting a way to correct it. The authors thank their institution for the support, under National Natural Science Foundation of China (No.~12071218), Young Elite Scientists Sponsorship Program by CAST (No.~156-O-110-0000085-3), Scientific Research Project of Colleges and Universities in Anhui Province (No. 2023AH050264), and Analysis, Probability \& Fractals Laboratory (No. LR18ES17).

\begin{spacing}{1}
	\bibliographystyle{amsplain}

\begin{thebibliography}{10}
	
	\bibitem{Anttila20223}  
R. Anttila, \textit{Pointwise Assouad dimension for measures}, Proc. Roy. Soc. Edinb. A \textbf{153} (2023), 2053--2078.

\bibitem{Anttila20224}  
R. Anttila, B. Bárány, and A. Käenmäki, \textit{Slices of the Takagi function}, Ergod. Theory Dyn. Syst. \textbf{44} (2024), 2361--2398.

\bibitem{Anttila2024}  
R. Anttila and A. Rutar, \textit{Fibre stability for dominated self-affine sets}, (2024), arXiv:2412.06579.

\bibitem{Assouad1983}  
P. Assouad, \textit{Plongements lipschitziens dans $\mathbb{R}^n$}, Bull. Soc. Math. Fr. \textbf{111} (1983), 429--448.

\bibitem{Banaji2024}  
A. Banaji and J. M. Fraser, \textit{Assouad type dimensions of infinitely generated self-conformal sets}, Nonlinearity \textbf{37} (2024), 045004.

\bibitem{Baranski2012} 
K. Baránski, \textit{On the dimension of graphs of Weierstrass-type functions with rapidly growing frequencies}, Nonlinearity \textbf{25} (2012), 193--209.

\bibitem{Baranski2014} 
K. Baránski,  B. Bárány, and J. Romanowska, \textit{On the dimension of the graph of the classical Weierstrass
	function}, Adv. Math. \textbf{265} (2014), 32--59.

\bibitem{Baranski2015}   K. Baránski, \textit{Dimension of the Graphs of the Weierstrass-Type Functions}, Cham: Springer (2015).

\bibitem{Besicovitch1937} 
A. S. Besicovitch and H. D. Ursell, \textit{Sets of fractional dimensions V: on dimensional numbers of some continuous curves}, J. Lond. Math. Soc. \textbf{32} (1937), 142--153.

\bibitem{Chen20222} 
H. P. Chen, J. M. Fraser, and H. Yu, \textit{Dimensions of the popcorn graph}, Proc. Amer. Math. Soc. \textbf{150} (2022), 4729--4742.

\bibitem{Chrontsios-Garitsis}  
E. K. Chrontsios-Garitsis and J. T. Tyson, \textit{On the Assouad spectrum of Hölder and Sobolev graphs}, Math. Proc. Cambridge Philos. Soc. \textbf{180} (2026), 105--131.

\bibitem{Falconer2014}  
K. J. Falconer, \textit{Fractal Geometry: Mathematical Foundations and Applications}, 3rd ed., Chichester: John Wiley Sons Ltd. (2014).

\bibitem{Feng2026} 
T. Y. Feng and J. M. Fraser, \textit{The Assouad spectrum and dimension of typical graphs}, (2026), arXiv:2601.05439.

\bibitem{Fraser2014}  
J. M. Fraser, \textit{Assouad type dimensions and homogeneity of fractals}, Trans. Amer. Math. Soc. \textbf{366} (2014), 6687--6733.

\bibitem{Fraser2021}  
J. M. Fraser, \textit{Assouad Dimension and Fractal Geometry}, Cambridge: Cambridge Univ. Press (2021).

\bibitem{Fraser2019}  
J. M. Fraser, K. E. Hare, K. G. Hare, S. Troscheit, and H. Yu, \textit{The Assouad spectrum and the quasi-Assouad dimension: a tale of two spectra}, Ann. Acad. Sci. Fenn. Math. \textbf{44} (2019), 379--387.

\bibitem{FraserYu2024}  
J. M. Fraser and A. Rutar, \textit{Assouad-type dimensions of overlapping self-affine sets}, Ann. Fenn. Math. \textbf{49} (2024), 3--21.

\bibitem{FraserYu2023}  
J. M. Fraser and L. Stuart, \textit{The Assouad spectrum of Kleinian limit sets and Patterson--Sullivan measure}, Geometriae Dedicata \textbf{217} (2023), 32.

\bibitem{FraserYu2018}  
J. M. Fraser and H. Yu, \textit{New dimension spectra: finer information on scaling and homogeneity}, Adv. Math. \textbf{329} (2018), 273--328.

\bibitem{Furstenberg2008} H. Furstenberg, \textit{Ergodic fractal measures and dimension conservation}, Ergod. Theory Dyn. Syst. \textbf{28} (2008), 405--422.

\bibitem{Hardy1916}  
G. H. Hardy, \textit{Weierstrass's non-differentiable function}, Trans. Amer. Math. Soc. \textbf{17} (1916), 301--325.

\bibitem{Howroyd2019}  
D. C. Howroyd and H. Yu, \textit{Assouad dimension of random processes}, Proc. Edinburgh Math. Soc. \textbf{62} (2019), 281--290.

\bibitem{Hunt1998}  
B. R. Hunt, \textit{The Hausdorff dimension of graphs of Weierstrass functions}, Proc. Amer. Math. Soc. \textbf{126} (1998), 791--800.

\bibitem{Hunt1991}  
B. R. Hunt and J. A. Yorke, \textit{Smooth dynamics on Weierstrass nowhere differentiable curves}, Trans. Amer. Math. Soc. \textbf{325} (1991), 141--154.

\bibitem{Jiang2025} 
L. Jiang, \textit{Assouad dimension of the Takagi function}, (2025), arXiv:2502.01140.

\bibitem{Kaplan1984} J. Kaplan, J. Mallet-Paret, and J. Yorke,  \textit{The Lyapunov dimension of a nowhere differentiable attracting
	torus}, Ergod. Theory Dyn. Syst. \textbf{4} (1984), 261--281.

\bibitem{Levy1934} P. L\'evy, \textit{Sur les int\'egrales dont les \'el\'ements sont des variables al\'eatoires ind\'ependentes}, Ann. Scuola. Norm.-Sci. \textbf{3} (1934), 337--366.

\bibitem{Lu2016} 
F. L\"u and L. Xi, \textit{Quasi-Assouad dimension of fractals}, J. Fractal Geom. \textbf{3} (2016),
187--215.

\bibitem{Ren2021}
H. J. Ren and W. X. Shen, \textit{A Dichotomy for the Weierstrass-type functions}, Invent. Math. \textbf{226} (2021), 1057--1100.

\bibitem{Shen2018}  
W. X. Shen, \textit{Hausdorff dimension of the graphs of the classical Weierstrass functions}, Math. Z. \textbf{289} (2018), 223--266.

\bibitem{Taylor1953}
 S. J. Taylor, \textit{The Hausdorff $ \alpha $-dimensional measure of Brownian paths in $ n $-space}, Math. Proc. Cambridge Philos. Soc. \textbf{49} (1953), 31--39.

\bibitem{Wang2025}  
H. Wang and J. Zahl, \textit{The Assouad dimension of Kakeya sets in $\mathbb{R}^3$}, Invent. Math. \textbf{241} (2025), 153--206.

\bibitem{Yu2024} 
B. Y. Yu and Y. S. Liang, \textit{On the Katugampola fractional integral and dimensional analysis of the fractal basin boundary for a random dynamical system}, Physica D \textbf{468} (2024), 134289.

\bibitem{Yu2020}  
H. Yu, \textit{Weak tangent and level sets of Takagi functions}, Monatsh. Math. \textbf{192} (2020), 249--264.
	
\end{thebibliography}

\end{spacing}
\end{document}